\newif\ifSIAM
\SIAMfalse

\ifSIAM
\documentclass[review,hidelinks,onefignum,onetabnum]{siamart220329}
\else
\documentclass[letterpaper]{article}
\usepackage[margin=1.2in]{geometry}
\usepackage{amsthm}
\fi

\usepackage{import}
\usepackage{mathtools}
\usepackage{amssymb}
\usepackage{mathrsfs}
\usepackage{xparse}
\usepackage{import}
\usepackage[normalem]{ulem}

\makeatletter
\@ifpackageloaded{unicode-math}{%
	\newcommand{\mymathbold}{\symbf}%
}{%
	\usepackage{bm}%
	\newcommand{\mymathbold}{\bm}%
}
\makeatother

\newcommand{\scrbar}[1]{\overline{\mathcal{#1}}}
\newcommand{\scrhat}[1]{\widehat{\mathcal{#1}}}
\newcommand{\scrtl}[1]{\widetilde{\mathcal{#1}}}

\ExplSyntaxOn

\cs_new_protected:Nn \bamboo_define:nnnnnN
{
	\cs_new_protected:cpx { #3 #1 #4 } { \exp_not:N #5{#6{#2}} }
}

\int_step_inline:nnn { `A } { `Z }
{
	\bamboo_define:nnnnnN
	{ \char_generate:nn { #1 } { 11 } }
	{ \char_generate:nn { #1 } { 11 } }
	{ bl }
	{ }
	{ \mymathbold }
	\use:n
	
	\bamboo_define:nnnnnN
	{ \char_generate:nn { #1 } { 11 } }
	{ \char_generate:nn { #1 } { 11 } }
	{ scr }
	{ }
	{ \mathcal }
	\use:n
	
	\bamboo_define:nnnnnN
	{ \char_generate:nn { #1 } { 11 } }
	{ \char_generate:nn { #1 } { 11 } }
	{ }
	{ hat }
	{ \widehat }
	\use:n
	
	\bamboo_define:nnnnnN
	{ \char_generate:nn { #1 } { 11 } }
	{ \char_generate:nn { #1 } { 11 } }
	{ scr }
	{ hat }
	{ \scrhat }
	\use:n
	
	\bamboo_define:nnnnnN
	{ \char_generate:nn { #1 } { 11 } }
	{ \char_generate:nn { #1 } { 11 } }
	{ }
	{ br }
	{ \overline }
	\use:n
	
	\bamboo_define:nnnnnN
	{ \char_generate:nn { #1 } { 11 } }
	{ \char_generate:nn { #1 } { 11 } }
	{ }
	{ tl }
	{ \widetilde }
	\use:n
	
	\bamboo_define:nnnnnN
	{ \char_generate:nn { #1 } { 11 } }
	{ \char_generate:nn { #1 } { 11 } }
	{ scr }
	{ br }
	{ \scrbar }
	\use:n
	
	\bamboo_define:nnnnnN
	{ \char_generate:nn { #1 } { 11 } }
	{ \char_generate:nn { #1 } { 11 } }
	{ scr }
	{ tl }
	{ \scrtl }
	\use:n
	
	\bamboo_define:nnnnnN
	{ \char_generate:nn { #1 } { 11 } }
	{ \char_generate:nn { #1 } { 11 } }
	{ }
	{ dt }
	{ \dot}
	\use:n
}
\int_step_inline:nnn { `a } { `z }
{
	\bamboo_define:nnnnnN
	{ \char_generate:nn { #1 } { 11 } }
	{ \char_generate:nn { #1 } { 11 } }
	{ bl }
	{ }
	{ \mymathbold }
	\use:n
	
	\bamboo_define:nnnnnN
	{ \char_generate:nn { #1 } { 11 } }
	{ \char_generate:nn { #1 } { 11 } }
	{ }
	{ hat }
	{ \hat }
	\use:n
	
	\bamboo_define:nnnnnN
	{ \char_generate:nn { #1 } { 11 } }
	{ \char_generate:nn { #1 } { 11 } }
	{ }
	{ br }
	{ \bar }
	\use:n
	
	\bamboo_define:nnnnnN
	{ \char_generate:nn { #1 } { 11 } }
	{ \char_generate:nn { #1 } { 11 } }
	{ }
	{ tl }
	{ \tilde }
	\use:n                            
	
	\bamboo_define:nnnnnN
	{ \char_generate:nn { #1 } { 11 } }
	{ \char_generate:nn { #1 } { 11 } }
	{ }
	{ dt }
	{ \dot}
	\use:n
}
\clist_map_inline:nn
{
	Gamma,Delta,Theta,Lambda,Xi,Pi,Sigma,Phi,Psi,Omega
}
{
	\bamboo_define:nnnnnN
	{ #1 }
	{ #1 }
	{ bl }
	{ }
	{ \mymathbold }
	\use:c
	
	\bamboo_define:nnnnnN
	{ #1 }
	{ #1 }
	{ op }
	{ }
	{ \op }
	\use:c
	
	\bamboo_define:nnnnnN
	{ #1 }
	{ #1 }
	{ scr }
	{ }
	{ \mathcal }
	\use:c
	
	\bamboo_define:nnnnnN
	{ #1 }
	{ #1 }
	{ }
	{ hat }
	{ \widehat }
	\use:c
	
	\bamboo_define:nnnnnN
	{ #1 }
	{ #1 }
	{ }
	{ br }
	{ \overline }
	\use:c
	
	\bamboo_define:nnnnnN
	{ #1 }
	{ #1 }
	{ }
	{ tl }
	{ \widetilde }
	\use:c
	
	\bamboo_define:nnnnnN
	{ #1 }
	{ #1 }
	{ }
	{ dt }
	{ \dot}
	\use:c
	
}
\clist_map_inline:nn
{
	alpha,beta,gamma,delta,epsilon,zeta,eta,theta,iota,kappa,
	lambda,mu,nu,xi,pi,rho,sigma,tau,phi,chi,psi,omega
}
{
	\bamboo_define:nnnnnN
	{ #1 }
	{ #1 }
	{ bl }
	{ }
	{ \mymathbold }
	\use:c
	
	\bamboo_define:nnnnnN
	{ #1 }
	{ #1 }
	{ }
	{ hat }
	{ \hat }
	\use:c
	
	\bamboo_define:nnnnnN
	{ #1 }
	{ #1 }
	{ }
	{ br }
	{ \bar }
	\use:c
	
	\bamboo_define:nnnnnN
	{ #1 }
	{ #1 }
	{ }
	{ tl }
	{ \tilde }
	\use:c
	
	\bamboo_define:nnnnnN
	{ #1 }
	{ #1 }
	{ }
	{ dt }
	{ \dot}
	\use:c
}

\ExplSyntaxOff

\DeclareMathOperator{\E}{\mathbf{E}}

\newcommand{\rank}{\operatorname{rank}}

\newcommand{\tr}{\operatorname{tr}}

\newcommand{\diag}{\operatorname{diag}}

\DeclarePairedDelimiter{\abs}{\lvert}{\rvert}

\DeclarePairedDelimiter{\parens}{(}{)}
\DeclarePairedDelimiter{\brackets}{[}{]}
\DeclarePairedDelimiterX{\ip}[2]{\langle}{\rangle}{#1,#2}

\DeclarePairedDelimiterXPP{\normsub}[2]{}{\lVert}{\rVert}{_{#2}}{#1}
\DeclarePairedDelimiterXPP{\ipsub}[3]{}{\langle}{\rangle}{_{#3}}{#1,#2}

\DeclarePairedDelimiterXPP{\ipHS}[2]{}{\langle}{\rangle}{_{\mathrm{HS}}}{#1, #2}
\DeclarePairedDelimiterXPP{\normHS}[1]{}{\lVert}{\rVert}{_{\mathrm{HS}}}{#1}

\DeclarePairedDelimiterXPP{\ipF}[2]{}{\langle}{\rangle}{_{\mathrm{F}}}{#1, #2}
\DeclarePairedDelimiterXPP{\normF}[1]{}{\lVert}{\rVert}{_{\mathrm{F}}}{#1}

\DeclarePairedDelimiterXPP{\dkl}[2]{\operatorname{D_{KL}}}{(}{)}{}{#1 \: \delimsize\Vert \: #2}

\DeclarePairedDelimiterXPP{\restr}[2]{}{{}}{\vert}{_{#2}}{#1}

\newcommand{\R}{\mathbf{R}}
\newcommand{\C}{\mathbf{C}}
\newcommand{\Z}{\mathbf{Z}}

\newcommand{\indicator}[1]{\mathbf{1}_{\{ #1 \}}}
\newcommand{\range}{\operatorname{range}}

\newcommand{\nul}{\operatorname{null}}

\newcommand{\normaldist}{\operatorname{\mathcal{N}}}

\newcommand{\negqquad}{\mspace{-36mu}}

\usepackage{array}
\usepackage{booktabs}
\usepackage{suffix}
\usepackage{bookmark}
\usepackage{color}
\usepackage{graphicx}
\graphicspath{{./figures/}}
\usepackage{subcaption}
\usepackage[capitalize]{cleveref}
\usepackage{ulem}
\usepackage{tikz}
\usepackage{pgfplots}
\pgfplotsset{compat=1.18}
\usetikzlibrary{arrows.meta}
\usepackage{hyperref}

\DeclarePairedDelimiterXPP{\opnorm}[1]{}{\lVert}{\rVert}{_{\ell_2}}{#1}
\DeclarePairedDelimiterXPP{\nucnorm}[1]{}{\lVert}{\rVert}{_{*}}{#1}

\newcommand{\SOgrp}{\mathrm{SO}}

\newcommand{\Ogrp}{\mathrm{O}}
\newcommand{\Or}{\Ogrp(r)}
\newcommand{\Ugrp}{\mathrm{U}}
\newcommand{\SUgrp}{\mathrm{SU}}
\newcommand{\SE}{\mathrm{SE}}
\newcommand{\St}{\mathrm{St}}
\newcommand{\st}{\text{ s.t.\ }}
\newcommand{\sbd}{\operatorname{SBD}}

\newcommand{\onevec}{\mathbf{1}}
\newcommand{\ER}{Erd\H{o}s--R\'enyi}
\newcommand{\real}{\operatorname{real}}

\ifSIAM
\title{Benign landscapes of low-dimensional relaxations%
	\\for orthogonal synchronization on general graphs%
	\thanks{Submitted to the editors July 6, 2023.
		\funding{This work was supported by the Swiss State Secretariat for Education, Research and Innovation (SERI) under contract number MB22.00027.}}}

\author{Andrew D.\ McRae%
	\thanks{Institute of Mathematics, EPFL, Lausanne (\email{andrew.mcrae@epfl.ch}, \email{nicolas.boumal@epfl.ch}).}%
	\and Nicolas Boumal\footnotemark[2]%
}

\headers{Benign landscapes for synchronization on graphs}{A. D. McRae and N. Boumal}

\else
\hypersetup{
	colorlinks,
	linkcolor=blue,
	citecolor=blue,
	urlcolor=blue,
}
\newtheorem{lemma}{Lemma}
\newtheorem{theorem}{Theorem}
\newtheorem{corollary}{Corollary}

\title{Benign landscapes of low-dimensional relaxations%
	\\for orthogonal synchronization on general graphs}

\author{Andrew D.\ McRae and Nicolas Boumal\thanks{Institute of Mathematics, EPFL, Lausanne (\texttt{andrew.mcrae@epfl.ch}, \texttt{nicolas.boumal@epfl.ch}). This work was supported by the Swiss State Secretariat for Education, Research and Innovation (SERI) under contract number MB22.00027.}}
\fi

\begin{document}
	
\maketitle

\begin{abstract}
Orthogonal group synchronization is the problem of estimating $n$ elements $Z_1, \ldots, Z_n$ from the $r \times r$ orthogonal group
given some relative measurements $R_{ij} \approx Z_i^{}Z_j^{-1}$.
The least-squares formulation is nonconvex.
To avoid its local minima, a Shor-type convex relaxation squares the dimension of the optimization problem from $O(n)$ to $O(n^2)$.
Alternatively, Burer--Monteiro-type nonconvex relaxations have generic landscape guarantees at dimension $O(n^{3/2})$.
For smaller relaxations, the problem structure matters.
It has been observed in the robotics literature that, for SLAM problems, it seems sufficient to increase the dimension by a small constant multiple over the original.
We partially explain this.
This also has implications for Kuramoto oscillators.

Specifically, we minimize the least-squares cost function in terms of estimators $Y_1, \ldots, Y_n$.
For $p \geq r$, each $Y_i$ is relaxed to the Stiefel manifold $\mathrm{St}(r, p)$ of $r \times p$ matrices with orthonormal rows.
The available measurements implicitly define a (connected) graph $G$ on $n$ vertices.
In the noiseless case, we show that, for all connected graphs $G$, second-order critical points are globally optimal as soon as $p \geq r+2$.
(This implies that Kuramoto oscillators on $\mathrm{St}(r, p)$ synchronize for all $p \geq r + 2$.)
This result is the best possible for general graphs; the previous best known result requires $2p \geq 3(r + 1)$.
For $p > r + 2$, our result is robust to modest amounts of noise (depending on $p$ and $G$).
Our proof uses a novel randomized choice of tangent direction to prove (near-)optimality of second-order critical points.
Finally, we partially extend our noiseless landscape results to the complex case (unitary group); we show that there are no spurious local minima when $2p \geq 3r$.
\end{abstract}

\section{Introduction and results}

We examine the optimization landscape of a class of quadratically constrained quadratic programs (QCQPs) that arise from the \emph{orthogonal group synchronization} problem.
This widely-studied problem has applications notably in simultaneous localization and mapping (SLAM) \cite{Carlone2015a}, cryo-electron microscopy (cryo-EM) \cite{Wang2013}, computer vision \cite{Martinec2007}, and phase retrieval \cite{Iwen2020}.
It also connects mathematically with oscillator networks \cite{Markdahl2020}.

Our main results in this paper are presented in this section as follows:
\Cref{sec:intro_group_est} presents the (real) orthogonal synchronization problem on a graph and our main optimization landscape results for the resulting QCQPs. \Cref{sec:intro_oscillators} connects our work to oscillator networks on (real) Stiefel manifolds and gives a new and optimal result in this field. \Cref{sec:intro_complex} partially extends these results to the complex case.
\Cref{sec:intro_smallr} details some specific implications for the low-dimensional groups that are of primary interest in many applications.
\Cref{sec:intro_addl} gives additional implications of our analysis that may be of independent interest.
\subsection{Problem setup and optimization landscape results}
\label{sec:intro_group_est}

The orthogonal synchronization problem we study is the following:
Let $G = (V, E)$ be a connected, undirected graph on the vertices $V = \{1, \dots, n\}$ for some integer $n \geq 1$.
Each vertex $i$ is associated with an unknown orthogonal matrix $Z_i \in \Or = \{ U \in \R^{r \times r} : UU^\top = I_r \}$.
We want to estimate $Z_1, \dots, Z_n$ from (potentially noisy) measurements of the form $R_{ij} = Z_i^{} Z_j^\top + \Delta_{ij} \in \R^{r \times r}$ for each edge $(i, j) \in E$,
where $\Delta_{ij}$ represents measurement error/noise.
Since the measurements are relative, estimation can only be done up to a global orthogonal transformation.

A simple least-squares\footnote{More generally, we could consider \emph{weighted} least-squares (e.g., if the noise variance changes across measurements). This corresponds to a graph $G$ with (positively) weighted edges.	Our analysis easily extends to this case, but we omit the generalization for simplicity.}
estimate of $Z_1, \dots, Z_n$ can be obtained from the following optimization problem:
\begin{equation}
	\label{eq:lsq_orig}
	\min_{Y \in \Or^n}~ \sum_{(i,j) \in E} \normF{Y_i - R_{ij} Y_j}^2,
\end{equation}
where $\normF{\cdot}$ denotes the matrix Frobenius (elementwise $\ell_2$) norm.
Although the cost function itself is convex in $Y$, the constraint set $\Or^n$ is nonconvex.
In general, the problem has spurious local minima in which local search methods (such as gradient descent) can get stuck.

Due to the orthogonality constraints, the above problem is equivalent to
\[
	\max_{Y \in \Or^n}~ \sum_{(i,j) \in E} \ip{R_{ij}}{Y_i^{} Y_j^\top},
\]
where $\ip{A}{B} = \tr(A B^\top)$ is the Hilbert--Schmidt (or Frobenius) matrix inner product.
To write this more compactly, let $C \in \R^{rn \times rn}$ denote the (incomplete) measurement matrix with blocks%
\footnote{Throughout this paper, indices into a matrix dimension of length $rn$ refer to blocks of $r$ rows or columns. Thus, for $C \in \R^{rn \times rn}$, $C_{ij}$ refers to the $(i,j)$th $r \times r$ block; for $Y \in \R^{rn \times p}$, $Y_i$ refers to the $i$th $r \times p$ block of $Y$, etc.}
\[
	C_{ij} = \begin{cases} R_{ij} & \textrm{ if } (i,j) \in E, \\ 0 & \text{ otherwise}. \end{cases}
\]
With the convention that $R_{ij} = R_{ji}^\top$, the matrix $C$ is symmetric.
Then, we can rewrite the problem as
\begin{equation}
	\label{eq:ncvx_orig}
	\max_{Y \in \R^{rn \times r}}~\ip {C}{Y Y^\top} \st Y_i^{} Y_i^\top = I_r, i = 1, \dots, n.
\end{equation}
The matrix $YY^\top$ is positive semidefinite with rank at most $r$.
Noting that $Y_i^{} Y_i^\top = (Y Y^\top)_{ii}$,
the classical (Shor) convex relaxation consists in replacing $YY^\top$ with a positive semidefinite matrix $X$, ignoring the rank constraint:
\begin{equation}
	\label{eq:sdp_orig}
	\max_{\substack{X \in \R^{rn \times rn} \\ X \succeq 0}}~\ip{C}{X} \st X_{ii} = I_r, i = 1,\dots, n.
\end{equation}
This allows the rank of $X$ to grow up to $rn$.
Alternatively, we can allow $Y$ to have $p \geq r$ columns in~\eqref{eq:ncvx_orig}.
This effectively allows $YY^\top$ to have rank up to $p$, providing a more gradual relaxation as we increase $p$.
This yields the rank-$p$ Burer--Monteiro relaxation:
\begin{equation}
	\label{eq:rankp_orig}
	\max_{Y \in \R^{rn \times p}}~\ip{C}{Y Y^\top} \st Y_i^{} Y_i^\top = I_r, i = 1, \dots, n.
\end{equation}
We can view the rank $p \geq r$ as a hyperparameter that interpolates the original nonconvex problem \eqref{eq:ncvx_orig} (in the case $p = r$) and the full SDP relaxation \eqref{eq:sdp_orig} (which corresponds to $p \geq rn$).

This paper primarily considers the \emph{optimization landscape} of the rank-relaxed nonconvex problem \eqref{eq:rankp_orig} for various values of $p$.
In particular, we ask:
\begin{itemize}
	\item[\textbf{Q:}] \emph{How large does $p$ need to be so \eqref{eq:rankp_orig} has no spurious local optima?}
\end{itemize}

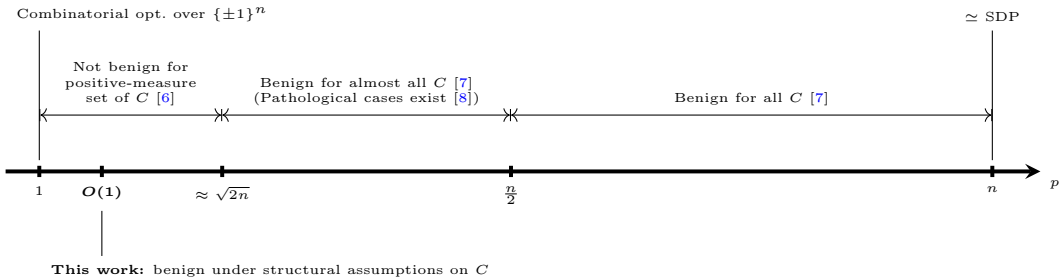
\begin{figure}[t]
	\centering
	\begin{tikzpicture}[font=\tiny]
		\begin{axis}[
			width=0.99\linewidth,
			height=60ex,
			xmin=-0.5, xmax=21,
			ymin=-3, ymax=7,
			axis x line=center,
			xlabel=$p$,%
			hide y axis,
			xtick={0.2, 1.5, 4, 10, 20},
			xticklabels={$1$, $\mymathbold{O(1)}$, $\approx \sqrt{2n}$, $\frac{n}{2}$, $n$},
			xticklabel style={align=center},
			xlabel style={below right},
			axis line style={line width=1.5pt}, 
			xtick style={line width=1.5pt, color=black} 
			]
			\draw[|<->|] (0.2, 1) -- node[above, align=center] {Not benign for\\ positive-measure \\ set of $C$ \cite{Waldspurger2020}} (4, 1);
			\draw[<->|] (4, 1) -- node[above, align=center] {Benign for almost all $C$ \cite{Boumal2019}\\ (Pathological cases exist \cite{OCarroll2022})} (10,1);
			\draw[<->|] (10, 1) -- node[above] {Benign for all $C$ \cite{Boumal2019}} (20,1);
			\draw (0.2, 0.2) -- (0.2, 2.5) node[above, align=center] {Comb\rlap{inatorial opt.\ over $\{\pm 1 \}^n$}};
			\draw (1.5, -0.7) -- (1.5, -1.5) node[below] {\textbf{This work:} \rlap{benign under structural assumptions on $C$}};
			\draw (20, 0.2) -- (20, 2.5) node[above] {$\simeq$ SDP};
		\end{axis}
	\end{tikzpicture}
	\caption{For low-dimensional relaxations, structural assumptions on the data (i.e., the cost matrix $C$) are necessary for benign nonconvexity. This diagram summarizes known results about the landscape of~\eqref{eq:rankp_orig} in the case $r = 1$, that is, relaxation at rank $p$ for synchronization of $n$ elements in $\Ogrp(1) = \{\pm 1\}$.}
	\label{fig:p_reqs}
\end{figure}

Much is known (see \Cref{fig:p_reqs} for a summary in the case $r = 1$):
For general cost matrices $C$, we need\footnote{Big-oh notation $O(\cdot)$ is with respect to $n$, treating $r$ as a small constant.} $p = O(n)$ to guarantee such a benign landscape (resulting in $O(n^2)$ variables, which is the same order as the full SDP relaxation);
with the assumption that $C$ is ``generic'' (outside of a zero measure set), we can reduce this to $p = O(n^{1/2})$ (resulting in $O(n^{3/2})$ variables);
for lower values of $p$, however, ``bad'' matrices $C$ are plenty: benign landscapes \emph{require} a structured $C$.

Remarkably, for cost matrices $C$ arising in specific applications, it has been observed that $p$ can be taken much smaller---just slightly above $r$.
This was observed and partially explained theoretically for models related to ours (complete-graph synchronization and stochastic block models) in prior works such as \cite{Bandeira2016a,Ling2022b}.
In the case of synchronization on a general graph for SLAM (robotics), similar empirical observations were reported in \cite{Rosen2019,Dellaert2020}; this is the inspiration for our work.

We prove that benign nonconvexity occurs for small values of $p$ for cost matrices $C$ arising from the general-graph synchronization problem; this gives additional theoretical support to the algorithms and empirical observations of \cite{Rosen2019,Dellaert2020}.
Setting $p = O(1)$ results in an optimization problem with $O(n)$ variables, similar to the original problem~\eqref{eq:ncvx_orig} but tractable despite being nonconvex.

Specifically, we give conditions under which every 
\emph{second-order critical point} (in particular, every local optimum) of~\eqref{eq:rankp_orig} is a global optimum.
We define such a point more precisely in \Cref{sec:math_background}, but it is, essentially, a point $Y \in \R^{rn \times p}$ where, subject to the constraints, the gradient at $Y$ is zero and all eigenvalues of the Hessian are nonpositive.

We first consider the case with no measurement error (i.e., $R_{ij} = Z_i^{} Z_j^\top$).
Clearly, a globally optimal solution to the original least-squares problem \eqref{eq:lsq_orig} is
\[
Z \coloneqq \begin{bmatrix*}
	Z_1 \\ \vdots \\ Z_n
\end{bmatrix*} \in \R^{rn \times r}.
\]
Absent noise, the computational task is trivial: fix $Z_1 = I_r$ arbitrarily (since recovery is up to global orthogonal transformation), then traverse any spanning tree of $G$,  recursively applying the measured relative differences to infer the other $Z_i$.
However, in general, the \emph{landscapes} of \eqref{eq:lsq_orig}, \eqref{eq:ncvx_orig} have spurious local optima.
What is the effect of relaxation then?

Our first main result states that, without noise, relaxing to $p = r+2$ eliminates spurious optima.
\begin{theorem}
	\label{thm:noiseless}
    Suppose $G$ is connected.
	If the measurements are exact, i.e., $R_{ij} = Z_i^{} Z_j^\top$ for all $(i,j) \in E$,
	then, if $p \geq r + 2$, any second-order critical point $Y$ of \eqref{eq:rankp_orig} satisfies $Y Y^\top = Z Z^\top$.
	Equivalently, $Y = Z U$ for some $r \times p$ matrix $U$ satisfying $U U^\top = I_r$.
\end{theorem}
This result is tight.
Indeed, at $p = r$ and $p = r+1$, there exist spurious local optima for certain connected graphs $G$;
see \Cref{sec:intro_oscillators} for further details and discussion.

We next consider the effect of measurement error (noise).
Our results do not apply to arbitrary noise (otherwise they would apply to arbitrary $C$)
but depend on the noise level relative to the graph connectivity.

To quantify the noise, denote by $\Delta$ the $rn \times rn$ matrix with the errors $\Delta_{ij}$ in the appropriate places for $(i,j) \in E$, setting $\Delta_{ij} = 0$ for $(i,j) \notin E$.
Thus, $\Delta$ is the portion of the cost matrix $C$ that is due to measurement error.
Our results are stated in terms of $\opnorm{\Delta}$ (where $\opnorm{\cdot}$ denotes matrix operator norm).

We quantify the connectivity of $G$ in a spectral sense.
Let $L = L(G)$ be the unnormalized graph Laplacian of $G$, defined as $L = \diag(A \onevec) - A$, where $A$ is the adjacency matrix of $G$ and $\onevec \in \R^n$ is the all-ones vector.
The matrix $L$ is positive semidefinite with eigenvalues $0 = \lambda_1 \leq \lambda_2 \leq \cdots \leq \lambda_n = \opnorm{L}$.
If $G$ is connected, $\lambda_2 > 0$.
The better $G$ is connected, the larger $\lambda_2$ is;
$\lambda_2$ is often called the \emph{algebraic connectivity} or \emph{Fiedler value} of $G$.

With these definitions in place, we can state our first noisy landscape result:
\begin{theorem}
	\label{thm:noisy_rankbd}
    Suppose $G$ is connected and $p > r + 2$.
	Define
	\[
    	C_p \coloneqq \frac{2(p + r - 2)}{p - r - 2}.
	\]
	Then any second-order critical point $Y$ of \eqref{eq:rankp_orig} satisfies
	\[
	\rank(Y) \leq r + 5 C_p^2 \parens*{ \frac{\opnorm{\Delta}}{\lambda_2} }^2 rn.
	\]
	If, furthermore, $p > r + 5 C_p^2 \parens*{ \frac{\opnorm{\Delta}}{\lambda_2} }^2 rn$,
	then $Y$ is a globally optimal solution to \eqref{eq:rankp_orig},
	and $Y Y^\top$ is an optimal solution to the SDP \eqref{eq:sdp_orig}.
\end{theorem}
This result quantitatively bounds how large $p$ needs to be so that \eqref{eq:rankp_orig} has a benign landscape and yields an exact solution to the full SDP relaxation.
The bound depends on the effective signal-to-noise ratio $\opnorm{\Delta} / \lambda_2$.

When \Cref{thm:noisy_rankbd} ensures that $\rank(Y) = r$ exactly, we obtain a stronger result:
\begin{theorem}
	\label{thm:noisy_landscape}
	Suppose $G$ is connected and $p > r + 2$.
	Let $C_p$ be as defined in \Cref{thm:noisy_rankbd}.
	If
	\begin{equation}
		\label{eq:noisy_landscape_cond}
		\opnorm{\Delta} < \frac{\lambda_2}{\sqrt{5} C_p \sqrt{rn}},
	\end{equation}
	then any second-order critical point $Y$ of \eqref{eq:rankp_orig} satisfies the following:
	\begin{itemize}
		\item $Y$ is the unique solution to \eqref{eq:rankp_orig} up to a global orthogonal transformation.
		\item $Y$ has rank $r$ and hence can be factored as $Y = \Zhat U$ for some $\Zhat \in \R^{rn \times r}$ and $U \in \R^{r \times p}$ such that $U U^\top = I_r$. Moreover, $\Zhat$ is the unique solution of \eqref{eq:ncvx_orig} up to a global orthogonal transformation.
		\item $Y Y^\top$ is the unique solution to the SDP \eqref{eq:sdp_orig}.
	\end{itemize}
\end{theorem}

First, note that if the measurement error $\Delta$ is small enough, then $p = r + 3$ suffices to obtain a benign landscape.
This provides baseline robustness for \Cref{thm:noiseless}.\footnote{The noiseless result \Cref{thm:noiseless} is already robust to noise when $p = r + 2$, but the required bound on $\opnorm{\Delta}$ might be much worse. See \Cref{sec:proof_noiseless} for details.}
Furthermore, the second-order critical points yield a global solution to the original, unrelaxed problem \eqref{eq:ncvx_orig}, which does not (in general) have a benign landscape.

If the measurement graph $G$ is complete, $\lambda_2 = n$, so the condition \eqref{eq:noisy_landscape_cond} becomes%
\footnote{Here and throughout the paper, $a \lesssim b$ ($a \gtrsim b$) means $a \leq c b$ ($a \geq c b$) for some unspecified constant $c > 0$. With subscripts, $a \lesssim_{p,r} b$ means $a \leq c(p, r) b$ (i.e., the ``constant'' could depend on $r$ and $p$).}
$\opnorm{\Delta} \lesssim_{p,r} \sqrt{n}$.
If the elements of $\Delta$ are, for example, i.i.d.\ Gaussian with variance $\sigma^2$,
this requires $\sigma^2 \lesssim_{p,r} 1$.
More generally, if $G$ is an \ER{} graph with edge probability $q$,
and $q \gtrsim \frac{\log n}{n}$ (which is necessary for $G$ to be connected),
we will have, with high probability, $\lambda_2 \approx nq$ (see \cite{Abdalla2022}).
In the i.i.d.\ Gaussian noise case, we have, with high probability, $\opnorm{\Delta} \approx \sigma \sqrt{nq}$ (see, e.g., \cite{Bandeira2016}),
so the condition on the noise variance becomes $\sigma^2 \lesssim_{p,r} q$.

In some cases, the assumptions of \Cref{thm:noisy_rankbd,thm:noisy_landscape} appear to be stronger than necessary (i.e., too pessimistic).
Our numerical experiments in \Cref{sec:sims} strongly suggest this for i.i.d.\ Gaussian noise.
It is also pessimistic compared to existing theoretical results for the complete-graph case.
Existing theory (e.g., \cite{Ling2022b}---see \Cref{sec:relwork_relax} for more references), under an additional assumption that the rows of $\Delta$ are not too correlated with $Z$, only requires $\opnorm{\Delta} \lesssim n^{3/4}$.
However, the techniques used are specialized to the complete-graph case,
and it is not clear how to adapt our methods in a way that recovers the best existing results when $G$ is the complete graph (see the end of \Cref{sec:proof_noisy_rankbound}).
Furthermore, it is unclear whether our requirements in \Cref{thm:noisy_rankbd,thm:noisy_landscape}
are still too pessimistic in the general adversarial-noise case.

\subsection{Implications for oscillator network synchronization}
\label{sec:intro_oscillators}
Another way to view the problem we have just described is \emph{oscillator synchronization}.
We briefly describe this connection here and spell out a corollary from \Cref{thm:noiseless}.
See, for example, \cite{Markdahl2018a,Markdahl2020} for more detailed discussion and derivations.

Given a connected graph\footnote{Again, the problem and our analysis easily extend to (positively) weighted edges.} $G$ defined as before,
a simple version of the \emph{Kuramoto} model for an oscillator network on $G$ is the following:
we have time-varying angles $\theta_1(t), \dots, \theta_n(t)$ associated with the $n$ vertices,
and these angles follow
\begin{equation}
	\label{eq:kuramoto_ang_eq}
	\dot{\theta}_i = -\sum_{j \in \scrN_i} \sin(\theta_i - \theta_j), \qquad i = 1, \dots, n,
\end{equation}
where $\scrN_i$ is the set of neighbors of vertex $i$ in $G$.
One can easily check that these dynamics are the gradient flow for the following optimization problem:
\begin{equation}
	\label{eq:kuramoto_ang_obj}
	\max_{\theta \in \R^n}~\frac{1}{2} \sum_{i,j = 1}^n A_{ij} \cos(\theta_i - \theta_j).
\end{equation}
Clearly, the optima of \eqref{eq:kuramoto_ang_obj} are the ``synchronized'' states $\theta_1 = \dots = \theta_n \mod 2\pi$.
Many papers have studied this model,
particularly with the following question in mind:
\begin{itemize}
	\item[\textbf{Q:}]
    \emph{For which graphs $G$ does the dynamical system \eqref{eq:kuramoto_ang_eq} converge to a synchronized state as $t \to \infty$ for ``generic'' initial conditions?}
\end{itemize}
Generic means ``except for a zero-measure set'' so that this happens with probability 1 if $\theta_1(0), \dots, \theta_n(0)$ are chosen uniformly at random over $[0, 2\pi)$.

To connect this problem to our work,
note that \eqref{eq:kuramoto_ang_obj} is a reparametrization of the problem \eqref{eq:rankp_orig} in the case $r = 1, p = 2$ when $Z_i = 1$ for all $i$ and the measurements are exact.\footnote{Explicitly, $\theta_i \mapsto Y_i = [\cos(\theta_i), \sin(\theta_i)]$ so that $\cos(\theta_i - \theta_j) = Y_i^{} Y_j^\top$. In this case, the cost matrix $C$ is equal to the adjacency matrix $A$ of $G$. The differential of the change of variable is surjective, hence it does not introduce spurious critical points \cite[Proposition 9.6]{Boumal2023}; thus the landscapes in $\theta$ and in $Y$ are qualitatively the same.}
There are many examples (see, e.g., \cite{Townsend2020}) of connected graphs $G$ under which \eqref{eq:kuramoto_ang_obj} has spurious local optima with strictly\footnote{Modulo the trivial direction that shifts all angles equally.} negative definite Hessian (and to which, therefore, gradient flow will converge if initialized in some positive-measure neighborhood).
This implies that \eqref{eq:rankp_orig} \emph{does not} always have a benign landscape for $r = 1, p = 2$.
However, this changes when one studies the ``synchronization'' landscape for higher-dimensional ``oscillators.''

For $r \geq 1$, consider \eqref{eq:rankp_orig} in the case that $Z_1 = \cdots = Z_n = I_r$.\footnote{This is without loss of generality as we can always smoothly change variables to bring the ground truth to this position without affecting the landscape: see \Cref{sec:real_proofs} for details.}
Furthermore, assume that there is no measurement error,\footnote{The simple connection to oscillators outlined here is most meaningful in the noiseless case.} so $R_{ij} = Z_i^{} Z_j^\top = I_r$ for all $(i,j) \in E$.
Note, furthermore, that the feasible points $Y = [Y_i]_i$ lie in a product of \emph{Stiefel manifolds} (we develop this connection further in the proofs of our main results):
\[
	\St(r, p) = \{ U \in \R^{r \times p} : U U^\top = I_r \}.
\]
With these simplifications and notation, problem~\eqref{eq:rankp_orig} becomes (within a factor of 2)
\begin{equation}
	\label{eq:kuramoto_stiefel_obj}
	\max_{Y \in \St(r,p)^n}~\frac{1}{2}\sum_{i,j=1}^n A_{ij} \tr(Y_i^{} Y_j^\top),
\end{equation}
where $A$ is the adjacency matrix of $G$.
Again, the global optima are precisely the $Y$ such that $Y_1 = \cdots = Y_n \in \St(r, p)$.

Any trajectory $Y(t)$ of the gradient flow on the constraint manifold satisfies the following system of differential equations:
\begin{equation}
	\label{eq:kuramoto_stiefel_eq}
	\Ydt_i = -\scrP_{T_{Y_i}} \parens*{ \sum_{j \in \scrN_i} (Y_i - Y_j) }, \qquad i = 1, \dots, n.
\end{equation}
Here, $T_U$ denotes the tangent space of $\St(r,p)$ at $U$ (see \Cref{sec:real_proofs}),
and $\scrP_{T_U}$ is the (Euclidean) orthogonal projection onto $T_U$.
These dynamics are (a simple version of) the Kuramoto model for a network of Stiefel-manifold valued oscillators.

Once again, we ask the question: when does the Stiefel-manifold valued Kuramoto oscillator network governed by \eqref{eq:kuramoto_stiefel_eq} \emph{synchronize} (i.e., converge to a synchronized state $Y_1 = \dots = Y_n$ as $t \to \infty$ for generic initialization)?
The remarkable results of \cite{Markdahl2018,Markdahl2018a,Markdahl2020} show that for certain manifolds, Kuramoto oscillator networks synchronize for \emph{any} connected graph $G$.
Specifically, the paper \cite{Markdahl2018} shows that for all $p \geq 3$, oscillator networks on $\St(1, p)$ (equivalently on the $d$-sphere $S^d$ for $d \geq 2$) always synchronize for connected $G$.
This was extended in \cite{Markdahl2018a,Markdahl2020} to show that $\St(r, p)$ oscillator networks synchronize if $2p \geq 3(r+1)$.

In this paper, we prove that $\St(r, p)$-valued oscillator networks synchronize under the weaker condition $p \geq r + 2$. This is a simple corollary of the noiseless orthogonal group synchronization result, \Cref{thm:noiseless}.
\begin{corollary}
	\label{cor:kuramoto_sync}
	For any connected graph $G$,
	the $\St(r, p)$-valued Kuramoto oscillator network on $G$ synchronizes if $p \geq r + 2$.
\end{corollary}
This result is tight and thus positively solves a conjecture in~\cite{Markdahl2020,Markdahl2021}.
More precisely, we have shown that general connected oscillator networks on $\St(r,p)$ synchronize if and only if the manifold is \emph{simply connected} (this is equivalent to $p \geq r + 2$: see, e.g., \cite[Example 4.53]{Hatcher2001}).
If the manifold is not simply connected, one can easily construct graphs $G$ such that \eqref{eq:kuramoto_stiefel_eq} has a stable equilibrium other than the synchronized state~\cite{Markdahl2021a}.

\Cref{cor:kuramoto_sync} follows from \Cref{thm:noiseless} by the same arguments used to prove \cite[Theorem 5.1]{Geshkovski2023}.
\Cref{thm:noiseless} implies that any critical point that is not the fully synchronized state is a strict saddle point of \eqref{eq:kuramoto_stiefel_obj} (i.e., the Riemannian Hessian has at least one strictly positive eigenvalue).
Because $\St(r, p)^n$ is compact, each gradient flow trajectory has at least one limit point (and any limit point is a critical point).
As $\St(r, p)^n$ and the objective function are real analytic, each trajectory must in fact converge to a single critical point (\cite{Lojasiewicz1983}, \cite[p.\ 4]{Lageman2007}).
By \cite[Lemma B.1]{Geshkovski2023},
the set of $Y(0)$ for which $Y(t)$ converges to a strict saddle has zero measure.
Then, for all $Y(0)$ outside this zero-measure set, $Y(t)$ must converge to a fully-synchronized state.

\subsection{The complex case}
\label{sec:intro_complex}
We can extend the previous (real) orthogonal matrix estimation problem to the complex case.
Here, we seek to estimate \emph{unitary} matrices $Z_1, \dots, Z_n \in \Ugrp(r) \coloneqq \{ U \in \C^{r \times r} : U U^* = I_r \}$
given measurements of the form $R_{ij} \approx Z_i^{} Z_j^*$.
We form our cost matrix $C \in \C^{rn \times rn}$ ($C$ is now \emph{Hermitian}, i.e., $C = C^*$)
and consider the relationships between the original unitary group least-squares problem
\begin{equation}
	\label{eq:cplx_ncvx}
	\max_{Y \in \C^{rn \times r}}~\ip {C}{Y Y^*} \st Y_i^{} Y_i^* = I_r, i = 1, \dots, n,
\end{equation}
the SDP relaxation
\begin{equation}
	\label{eq:cplx_sdp}
	\max_{\substack{X \in \C^{rn \times rn} \\ X \succeq 0}}~\ip{C}{X} \st X_{ii} = I_r, i = 1,\dots, n,
\end{equation}
and the rank-$p$ relaxation
\begin{equation}
	\label{eq:cplx_rankp}
	\max_{Y \in \C^{rn \times p}}~\ip {C}{Y Y^*} \st Y_i^{} Y_i^* = I_r, i = 1, \dots, n.
\end{equation}
We denote the complex $r \times p$ Stiefel manifold by
\[
	\St(r, p, \C) = \{ U \in \C^{r \times p} : U U^* = I_r \}.
\]
For simplicity, we only consider the noiseless landscape.
Again, the result has implications for Kuramoto oscillators on the complex Stiefel manifold,
which are related to the ``quantum Kuramoto model'' on unitary matrices (\cite{Lohe2009}; see \Cref{sec:relwork_oscillators} for additional references).
\begin{theorem}
	\label{thm:complex}
	Suppose $G$ is connected, and the measurements are exact (i.e., $R_{ij} = Z_i^{} Z_j^*$ for $(i,j) \in E$).
	If $2p \geq 3r$,	
    then any second-order critical point $Y$ of \eqref{eq:cplx_rankp} satisfies $Y Y^* = Z Z^*$.
	Consequently, the $\St(r, p, \C)$-valued Kuramoto oscillator network on $G$ synchronizes.
\end{theorem}

Due to the $3/2$ factor and certain similarities in the proof, this can also be seen as a complex adaptation of the result in \cite{Markdahl2020}.
Curiously, the innovations that allow us to improve that previous result in the real case do not easily carry over to the complex case; see the proof in \Cref{sec:complexcase} for further details.

As discussed in \Cref{sec:intro_oscillators}, our results for the real case show that we obtain a benign landscape (and all connected oscillator networks synchronize) as soon as the real Stiefel manifold $\St(r, p)$ is \emph{simply connected}.
In the complex case, the Stiefel manifold $\St(r, p, \C)$ is simply connected as soon as $p \geq r + 1$ (again, see \cite[Example 4.53]{Hatcher2001}).
We conjecture that this is the correct condition and that therefore \Cref{thm:complex} is suboptimal.
However, we have been unable to prove this, and it is unclear how to test it empirically.%
\footnote{The standard counterexamples to benign landscape/synchronization results are cycle graphs (indeed, the paper \cite{Markdahl2021a} uses cycle graphs to show that networks of oscillators taking values in a non--simply-connected manifold do not synchronize in general). We do not expect this counterexample to work in the complex case when $p \geq r + 1$, because a cycle graph corresponds geometrically to the unit circle, that is, $\Ugrp(1)$, and our result shows that, if $r = 1$, $p \geq 2 = r + 1$ suffices in the complex case. It is not clear how to construct a higher-dimensional equivalent as a candidate counterexample to our conjecture. For example, intuition from homotopy theory would suggest that we try a graph corresponding geometrically to a higher-dimensional sphere.}
If our conjecture is true, it is not clear where the limitations are in our current proof.
The critical step in our and others' proofs is to make an (educated) guess of potential descent direction for the objective function; how to improve this guess is unclear.

\subsection{Discussion of small-$r$ synchronization conditions}
\label{sec:intro_smallr}
It is interesting to consider the implications of the noiseless landscape results \Cref{thm:noiseless,cor:kuramoto_sync,thm:complex}
for small values of the matrix dimension $r$, as this covers many applications.
See \Cref{tab:groups} for a summary.

For synchronization of \emph{rotations} (rather than orthogonal transformations) in the plane $\R^2$, we can adopt two perspectives.
We could view $\SOgrp(2)$ as one of the two connected components of $\Ogrp(2)$, in which case \Cref{thm:noiseless} allows us to relax to synchronization on $\St(4, 2)$ (a 5-dimensional manifold).
Alternatively, we can view $\SOgrp(2)$ as isomorphic to $\Ugrp(1)$ (a circle, $S^1$), in which case \Cref{thm:complex} allows us to relax to synchronization on the \emph{complex} Stiefel manifold $\St(2, 1, \C)$, which is isomorphic to $S^3$ (a 3-dimensional sphere).

We have a similar situation for synchronization of rotations in $\R^3$.
Viewing $\SOgrp(3)$ as a subgroup of $\Ogrp(3)$,
we can relax the synchronization to $\St(5, 3)$, which is a 9-dimensional manifold.
Alternatively, it is well known that we can embed\footnote{More precisely, there is a double covering group homomorphism from $\SUgrp(2)$ to $\SOgrp(3)$.} $\SOgrp(3)$ in the special unitary group $\SUgrp(2)$,
which is the subgroup of $\Ugrp(2)$ of matrices with determinant 1.
With this formulation, we can synchronize in $\Ugrp(2)$ by relaxing to $\St(3, 2, \C)$, which is an 8-dimensional manifold.

\begin{table}[t]
	\centering
	\begin{tabular}{cccccc}
		\toprule
		Group & Field & $r$ & Min.\ $p$ req. & Orig. dim. & Min.\ Stiefel dim. \\
		\midrule
		$\Ogrp(1) = \Z_2$ & $\R$ & $1$ & $3$ & $0$ & $2$ \\
		$\Ugrp(1) = S^1 = \SOgrp(2)$ & $\C$ & 1 & 2 & 1 & 3 \\
		$\Ogrp(2)$ & $\R$ & $2$ & $4$ & $1$ & 5 \\
		$\Ugrp(2)$ & $\C$ & $2$ & $3$ & $4$ & $8$ \\
		$\Ogrp(3)$ & $\R$ & $3$ & $5$ & $3$ & $9$ \\
		$\Ugrp(3)$ & $\C$ & $3$ & $5$ ($4^*$) & $9$ & $21$ ($15^*$) \\
		\bottomrule
	\end{tabular}
	\caption{Properties of the groups $\Ogrp(r)$ and $\Ugrp(r)$ for $r \leq 3$ and the relaxations required to guarantee synchronization. In the last row, $^*$ indicates the result of the conjectured $p \geq r + 1$ condition for the complex case.}
	\label{tab:groups}
\end{table}

\subsection{Additional results}
\label{sec:intro_addl}
Our analysis yields several additional results of independent interest.
For the nonconvex relaxation~\eqref{eq:rankp_orig}, we show second-order critical points can be good approximate solutions even if they are not globally optimal.
For the SDP relaxation~\eqref{eq:sdp_orig}, we show a tightness result (rank recovery) and an approximation result.

\subsubsection{Error bounds for all second-order critical points}
\label{sec:intro_errorbounds}
Regardless of whether the conditions of \Cref{thm:noisy_rankbd,thm:noisy_landscape} are satisfied,
we can obtain useful error bounds for all second-order critical points $Y$ of \eqref{eq:rankp_orig}.
Recall that $\opnorm{\cdot}$, $\normF{\cdot}$ respectively denote the operator and Frobenius matrix norms.
The quality metric we use for a candidate solution $Y$ is the correlation $\ip{Z Z^\top}{Y Y^\top} = \normF{Z^\top Y}^2$.
The maximum value this can take is $n^2 r$ (because $\opnorm{Z Z^\top} = n$, and $\tr(Y Y^\top) = rn$),
and this value is reached if and only if $Y Y^\top = Z Z^\top$ (as in \Cref{thm:noiseless}).
\begin{theorem}
	\label{thm:noisy_bound}
	Assume $G$ is connected and $p > r + 2$.
    Then any second-order critical point $Y$ of \eqref{eq:rankp_orig} satisfies
	\[
		\ip{Z Z^\top}{Y Y^\top} \geq \brackets*{ 1 - C_p^2 \frac{\opnorm{\Delta}^2}{\lambda_2^2} }n^2 r,
	\]
	where $C_p$ is defined in the statement of \Cref{thm:noisy_rankbd}.
\end{theorem}

Comparable error bounds have been shown for the eigenvector method;
in particular, our bound is identical (within constants and dependence on $p$) to one in \cite{Doherty2022}.
See \Cref{sec:relwork_spectral} for further references.
In the \ER{} graph case with i.i.d.\ zero-mean random noise (see discussion after \Cref{thm:noisy_landscape}),
we obtain, with high probability,
\[
	1 - \frac{\ip{Z Z^\top}{Y Y^\top}}{n^2 r} \lesssim \frac{\sigma^2}{q n},
\]
where $q$ is the edge probability and $\sigma^2$ is the noise variance.
This is the minimax-optimal error rate for this problem, as shown in \cite{Zhang2022} (within constants and dependence on $r$---the error metric in \cite{Zhang2022} is slightly different and depends differently on $r$). Again, see \Cref{sec:relwork_spectral} for more references.

\subsubsection{Consequences for the SDP relaxation}
For the SDP relaxation~\eqref{eq:sdp_orig}, we first have an exactness result (tight relaxation, rank recovery):
\begin{corollary}
	\label{cor:sdp_exactness}
    Assume $G$ is connected. If $\opnorm{\Delta} < \frac{\lambda_2}{2\sqrt{5} \sqrt{rn}}$,
	then
	\begin{itemize}
		\item The SDP relaxation \eqref{eq:sdp_orig} has a unique solution $\Xhat$, and $\rank(\Xhat) = r$.
		\item The unrelaxed problem \eqref{eq:ncvx_orig} has a unique (up to orthogonal transformation) solution $\Zhat \in \R^{rn \times r}$.
		\item $\Xhat = \Zhat \Zhat^\top$.
	\end{itemize}
\end{corollary}

Next, we have a general error bound for the SDP relaxation that applies even if the exactness result above does not:
\begin{corollary}
	\label{cor:sdp_bound}
    Assume $G$ is connected.
	Any solution $X$ to \eqref{eq:sdp_orig} satisfies
	\[
	\ip{X}{Z Z^\top} \geq \parens*{ 1 - \frac{4 \opnorm{\Delta}^2}{\lambda_2^2} } n^2 r.
	\]
\end{corollary}
These corollaries follow from \Cref{thm:noisy_landscape,thm:noisy_bound} in the limit $p \to \infty$.
To be precise, let $X$ be an optimal solution to \eqref{eq:sdp_orig}.
For any $p \geq rn$, there exists $Y \in \R^{rn \times p}$ such that $X = Y Y^\top$.
The fact that $X$ is feasible implies that $Y$ is feasible for \eqref{eq:rankp_orig}.
Furthermore, the optimality of $X$ implies that $Y$ is a global optimum and therefore is a second-order critical point \cite[Prop.~2.4]{Boumal2019}.
We then apply \Cref{thm:noisy_landscape,thm:noisy_bound} and take $p \to \infty$, noting that $C_p \to 2$.
See the discussion in \Cref{sec:intro_errorbounds,sec:relwork_spectral}
for comparison to existing error bounds.
\section{Related work}
\label{sec:relwork}

The literature on orthogonal group synchronization is vast, appearing in multiple communities such as robotics, image processing, signal processing, and dynamical systems.
We highlight a few salient references here; see also \cite{Hartley2013,Rosen2021} for partial surveys.
Many of the tools we use in our analysis have been used before.
We point this out in our analysis along the way.

\subsection{Rank relaxation for synchronization}
\label{sec:relwork_relax}
Low-rank factorizations of SDPs (which, in our case, correspond to partial rank relaxations of the synchronization problem)
have a long history.
This approach is often called Burer--Monteiro factorization after the pioneering work of those authors (e.g., \cite{Burer2002,Burer2005}).

The report \cite{Boumal2015}, along with the more general results in \cite{Boumal2019}, provides a theoretical framework for analyzing Burer--Monteiro factorizations (like ours) of SDPs with (block-)diagonal constraints. The papers \cite{Rosen2019} and \cite{Dellaert2020} develop fast algorithms (generalized to the special Euclidean group case in \cite{Rosen2019}), showing that the ``Riemannian staircase'' approach (iteratively increasing the relaxation rank) proposed in \cite{Boumal2015} provides an exact solution to the SDP relaxation, but they do not specify what relaxation rank ($p$ in our notation) suffices. Our results provide an upper bound (optimal in the noiseless case) on how much such algorithms must relax the rank constraint for the synchronization problem.
 
As far as we are aware, \emph{landscape} results similar to ours have previously only been proved in the complete measurement graph case.
For this case, the papers \cite{Bandeira2016a,Ling2022b} provide error bounds and benign landscape results for rank-relaxed optimization like \eqref{eq:rankp_orig} (for the case $r = 1$, $p = 2$ in \cite{Bandeira2016a}).
 
The paper \cite{Ling2022b} analyzes the landscape of \eqref{eq:rankp_orig} in the same general $\Or$ case that we do and is thus the most comparable work to ours.
We can directly compare our \Cref{thm:noisy_landscape} with \cite[Theorem 3]{Ling2022b}.
In the complete-graph case, $\lambda_2 = n - 1$, so the condition \eqref{eq:noisy_landscape_cond} of \Cref{thm:noisy_landscape} becomes $\opnorm{\Delta} \lesssim_{p,r} \sqrt{n}$.
This prior result \cite[Theorem 3]{Ling2022b}, in the adversarial-noise case (i.e., only assuming a bound on $\opnorm{\Delta}$),
has a comparable requirement.
With additional assumptions on $\Delta$, the prior result improves this to $\opnorm{\Delta} \lesssim n^{3/4}$,
but the techniques used do not easily carry over to our general-graph case.
In addition, the paper \cite{Ling2022b} requires $p > 2r$ for all its results.
 
To the best of our knowledge, the partial benign landscape result of \Cref{thm:noisy_rankbd} (which gives conditions under which \eqref{eq:rankp_orig} has a benign landscape and yields an exact solution to the SDP even when the solution is not necessarily rank-$r$) is new even in the complete-graph case.

In a different vein, \cite{Mei2017} provides a general bound on how well rank-$p$ Burer--Monteiro factorizations for $\Ogrp(r)$ optimization can approximate the full SDP relaxation in terms of \emph{objective function value}.
Their results bound the approximation error by a term proportional to $r/p$.
Our results show that, in some cases, we already obtain perfect approximation with $p$ only slightly larger than $r$.

An interesting parallel work on low-rank Burer--Monteiro factorization landscapes is \cite{Zhang2022a}.
However, their setting is quite different from ours (requiring a strongly convex objective and no constraints), so their results are not directly comparable.

\subsection{The spectral approach and previous error bounds}
\label{sec:relwork_spectral}
Our work follows a large body of prior results that use spectral properties of the measurement graph.
In particular, these results use the eigenvalues and eigenvectors of the graph Laplacian matrix $L$ (or, in some cases, the adjacency matrix $A$) of the graph $G$.
A particularly important quantity is the \emph{graph connection Laplacian} matrix (defined in \cite{Singer2012} for different purposes) which can be formed directly from our observations $R_{ij}$ (this is precisely the matrix $\Lhat$ in our notation---see \Cref{sec:math_background} for the definition).

The \emph{eigenvector method} was introduced in \cite{Singer2011} for the purpose of rotation synchronization;
this method directly uses the eigenvectors of the graph connection Laplacian (or, in the original paper, the adjacency matrix equivalent).
Closely related to this, the paper \cite{Bandeira2013} studies the relationship between the eigenvalues/vectors of the graph connection Laplacian and the optimal objective function value of \eqref{eq:ncvx_orig}.

The paper \cite{Bandeira2018}, by analyzing eigenvalues of $\Lhat$, provides conditions under which the SDP relaxation gives exact recovery for $\Z_2$ (i.e., $\Ogrp(1)$) synchronization on an \ER{} graph.
The paper \cite{Wang2013a} considers a robust version of the SDP relaxation that uses the sum of absolute errors rather than least-squares.
They provide conditions for exact recovery for an \ER{} random graph and sparse errors.

The papers \cite{Iwen2020,Filbir2021,Doherty2022} and \cite[Chapter 6]{Preskitt2018} contain a variety of error bounds for the eigenvector method in terms of the graph Fiedler value $\lambda_2(L)$ and various norms of the measurement error matrix $\Delta$. 
This is extended to $\SE(r)$ synchronization (special Euclidean group) in \cite{Doherty2022}.
These results are the most comparable to our error bounds, because they apply to general graphs and measurement errors.
The bound in \cite{Doherty2022} is, within constants, identical to our bounds in \Cref{thm:noisy_bound,cor:sdp_bound}.

The papers \cite{Gao2021,Gao2022,Zhang2022} show that the SDP relaxation, the (rounded) eigenvector method, and generalized power method (eigenvector method followed by iterative refinement) all achieve asymptotically minimax-optimal error with Gaussian noise on Erdős--Rényi random graphs (interestingly, \cite{Zhang2022}, like the older paper \cite{Singer2012}, analyzes the adjacency matrix leading eigenvector). These results agree with our error bounds within constants (see \Cref{sec:intro_errorbounds}), though our results apply to much more general situations.

\subsection{Oscillator synchronization literature}
\label{sec:relwork_oscillators}
Oscillator synchronization is a large field of research. See \cite{Ha2016,Pikovsky2023} for surveys. Our work touches the small subset corresponding to simplified Kuramoto oscillators.

As discussed in \Cref{sec:intro_oscillators}, the classical and most well-studied Kuramoto oscillator is that of angular synchronization (i.e., on the unit circle $S^1$).
One line of research studies \emph{which connected graphs} synchronize on $S^1$.
This includes deterministic guarantees based on the density of a graph \cite{Taylor2012,Ling2019,Townsend2020,Yoneda2021,Lu2020,Kassabov2021} and high-probability properties of random graphs \cite{Ling2019,Kassabov2022,Abdalla2022,Abdalla2022a}.
For example, \cite{Kassabov2021} shows that every graph synchronizes in which every vertex is connected to at least 3/4 of the other vertices.
The paper \cite{Abdalla2022} gives more general deterministic conditions based on spectral expander properties and shows that, asymptotically, \ER{} graphs that are connected also synchronize.

A complementary line of research studies \emph{for which manifolds} do \emph{all} connected oscillator networks synchronize.
It was shown in \cite{Markdahl2018,Markdahl2021} that Kuramoto networks on the $d$-sphere $S^d$ synchronize for any $d \geq 2$.
This corresponds to the Stiefel manifold $\St(r, p)$ for $r = 1$ and $p \geq 3$.
More generally, the papers \cite{Markdahl2018a,Markdahl2020}
show that networks on $\St(r, p)$ synchronize if $2p \geq 3(r + 1)$.
Our \Cref{cor:kuramoto_sync} improves this condition to $p \geq r + 2$,
which is optimal according to \cite{Markdahl2018b,Markdahl2021a}
and confirms a conjecture in \cite{Markdahl2020,Markdahl2021}.
See the discussion around \Cref{cor:kuramoto_sync} for further details.
See also \Cref{sec:intro_complex} for similar discussion in the complex case.

Oscillators on complex Stiefel manifolds can be seen as a generalization of the ``quantum'' Kuramoto oscillators introduced in \cite{Lohe2009}. See \cite{Ha2016,Deville2019} for more references.

\section{Key mathematical tools}
\label{sec:math_background}
In this section, we make precise and fill out the mathematical framework for our analysis and results.
We only consider the real case in this section and the next.
We make the necessary adjustments for the complex case in \cref{sec:complexcase}.

\subsection{Graph Laplacian formulation}
In each optimization problem we consider,
the orthogonality/block diagonal constraint ensures that the $r\times r$ diagonal blocks of the cost matrix $C$ have no effect.
We therefore replace $C$ by another matrix that will be more convenient for analysis.

Let $A \in \R^{n \times n}$ be the adjacency matrix of $G$ (i.e., the symmetric matrix with $A_{ij} = \indicator{(i,j) \in E}$).
Let $L$ be the graph Laplacian matrix defined by
\[
L_{ij} = \begin{cases} \sum_{k \neq i} A_{ik} & \textrm{ if } i = j, \\ -A_{ij} & \textrm{ if } i \neq j. \end{cases}
\]
It is well known that $L$ is a positive semidefinite (PSD) matrix whose smallest eigenvalue is $\lambda_1(L) = 0$ with corresponding eigenvector $v_1 = \mathbf{1}_n$. The measurement graph is connected if and only if the second-smallest eigenvalue $\lambda_2 = \lambda_2(L) > 0$.

Recall that $Z_1, \dots, Z_n \in \Ogrp(r)$ are the ground-truth matrices that we want to estimate.
Let
\[
D = D(Z) \coloneqq \begin{bmatrix*}
	Z_1 & & & \\ & Z_2 & & \\ & & \ddots & \\ & & & Z_n
\end{bmatrix*}.
\]
With $\otimes$ denoting Kronecker product, let
\[
    L_Z \coloneqq D (L \otimes I_r) D^\top \quad \text{and} \quad \Lhat \coloneq L_Z - \Delta,
    \label{eq:LZLhat}
\]
where $\Delta$ is the symmetric matrix containing the measurement noise blocks $\Delta_{ij}$ (with $\Delta_{ij} = 0$ if $(i,j) \notin E$).
Note that the eigenvalues of $L_Z$ are identical to those of $L$ (with multiplicities multiplied by $r$),
and, if $G$ is connected, the $r$-dimensional subspace corresponding to $\lambda_1(L) = 0$ is precisely the span of the columns of $Z$.
The matrix $\Lhat$ was introduced in \cite{Singer2012} under the name \emph{graph connection Laplacian}.

Note that $C_{ij} = -\Lhat_{ij}$ for $i \neq j$,
and the diagonal blocks of the cost matrix have no effect on the optimization landscape of problems \eqref{eq:sdp_orig} and \eqref{eq:rankp_orig} due to the block diagonal constraint.
Therefore, the SDP \eqref{eq:sdp_orig} has the exact same landscape in the variable $X$ as
\begin{equation}
	\label{eq:sdp_laplace_blockwise}
	\min_{X \succeq 0}~\ip{\Lhat}{X} \st X_{ii} = I_r, i = 1,\dots, n.
\end{equation}
Likewise, the (relaxed) nonconvex problem \eqref{eq:rankp_orig} has the same landscape as
\begin{equation}
	\label{eq:rankp_laplace_blockwise}
	\min_{Y \in \R^{rn \times p}}~\ip {\Lhat}{Y Y^\top} \st Y_i^{} Y_i^\top = I_r, i = 1, \dots, n.
\end{equation}
From now on, we use those formulations.

\subsection{Manifold of feasible points and necessary optimality conditions}
\label{sec:manifold_desc}
Optimization problems of the form \eqref{eq:rankp_laplace_blockwise} have been well studied.
See, for example, \cite{Boumal2015,Boumal2019} for an overview.
We summarize the relevant facts in this section.

First, note that the constraints in \eqref{eq:sdp_laplace_blockwise} and \eqref{eq:rankp_laplace_blockwise} apply to the $r \times r$ diagonal blocks of a matrix.
To simplify notation, we denote the symmetric block-diagonal projection $\sbd \colon \R^{rn \times rn} \to \R^{rn \times rn}$ by
\[
	\sbd(X)_{ij}
	= \begin{cases}
		\frac{X_{ii} + X_{ii}^\top}{2} & \textrm{ if } i = j, \\
		0 & \textrm{ if } i \neq j.
	\end{cases}
\]
We can then write the semidefinite problem \eqref{eq:sdp_laplace_blockwise} more compactly as
\begin{equation}
	\label{eq:sdp_laplace_sbd}
	\min_{X \succeq 0}~\ip{\Lhat}{X} \st \sbd(X) = I_{rn}.
\end{equation}
Similarly, we can write \eqref{eq:rankp_laplace_blockwise} as
\begin{equation}
	\label{eq:rankp_laplace_sbd}
	\min_{Y \in \R^{rn \times p}}~\ip{\Lhat}{Y Y^\top} \st \sbd(Y Y^\top) = I_{rn}.
\end{equation}
The \emph{symmetrizing} aspect of the projection operator $\sbd$ is, for now, completely redundant but will become useful as we continue.

The feasible points of \eqref{eq:rankp_laplace_sbd} form a smooth submanifold $\scrM$ of $\R^{rn \times p}$:
\[
    \scrM \coloneq \St(r, p)^n = \{ Y \in \R^{rn \times p} : Y_i^{} Y_i^\top = I_r \textrm{ for } i = 1, \ldots, n\}.
\]
The most important object for us to understand is the \emph{tangent space} $T_Y$ at a point $Y \in \scrM$.
An $rn \times p$ matrix $\Ydt$ is in $T_Y$ if and only if its $r \times p$ blocks satisfy
$\Ydt_i^{} Y_i^\top + Y_i^{} \Ydt_i^\top = 0 \in \R^{r \times r}$ for all $i = 1, \dots, n$.
The orthogonal projection of an arbitrary $W \in \R^{rn \times p}$ onto $T_Y$ is given by
\[
	\scrP_{T_Y}(W) = W - \sbd(W Y^\top) Y.
\]
Local minima of \eqref{eq:rankp_laplace_sbd} are also first- and second-order critical \cite[Def.\ 2.3]{Boumal2019}: these conditions encode the fact that the Riemannian gradient is zero and the Riemannian Hessian is positive semidefinite (on the tangent space).
To be specific, let
\begin{align}
    S(Y) \coloneqq \Lhat - \sbd(\Lhat Y Y^\top).
    \label{eq:SY}
\end{align}
Then, one can show (see \cite{Boumal2019}) that any local minimum $Y$ satisfies both of the following:
\begin{itemize}
	\item First order condition: $S(Y) Y = 0$. This is proportional to the $T_Y$-projected gradient of the objective function.
	\item Second order condition: for all $\Ydt \in T_Y$, $\ip{S(Y) \Ydt}{\Ydt} \geq 0$. This says that the Riemannian Hessian of the objective function is positive semidefinite.
\end{itemize}
A point $Y$ that satisfies the first condition is a \emph{first-order critical point}.
A point $Y$ that satisfies both conditions is a \emph{second-order critical point}.
Note that second-order criticality is independent of whether we write the problem as \eqref{eq:rankp_orig} or \eqref{eq:rankp_laplace_sbd}
(in fact, $S(Y)$ is unchanged if we replace $\Lhat$ by $-C$).

\subsection{Dual certificates and global optimality}
\label{sec:dualcerts}
A natural question that remains is the following: how do we know if a candidate solution $Y$ to \eqref{eq:rankp_laplace_sbd} is, in fact, globally optimal?
This is a well-studied problem (again, see \cite{Boumal2019} for further discussion and references), and here we list some facts that will be useful throughout our analysis.

To study optimality conditions for \eqref{eq:rankp_laplace_sbd}, it is useful to consider the full SDP relaxation \eqref{eq:sdp_laplace_sbd}.
To show that a feasible point $X$ of \eqref{eq:sdp_laplace_sbd} is optimal,
it suffices to provide a \emph{dual certificate}.
For \eqref{eq:sdp_laplace_sbd}, a dual certificate is a matrix $S$ of the form
\[
	S = \Lhat - \sbd(\Lambda)
\]
for some $rn \times rn$ matrix $\Lambda$ with the conditions that (a) $S \succeq 0$ and (b) $\ip{S}{X} = 0$ (or $S X = 0$, which is equivalent because $S,X \succeq 0$).
If such a matrix $S$ exists, we have, for any feasible point $X' \succeq 0$,
\begin{align*}
	\ip{\Lhat}{X'} - \ip{\Lhat}{X}
	&= \ip{\Lhat}{X' - X} \\
	&= \ip{S}{X' - X} + \ip{\sbd(\Lambda)}{X' - X} \\
	&= \ip{S}{X'} - \underbrace{\ip{S}{X}}_{= 0} + \ip{\Lambda}{\underbrace{\sbd(X' - X)}_{= 0}}
	= \ip{S}{X'}
	\geq 0.
\end{align*}
Thus     $X$ is optimal. Furthermore, if $X'$ is also optimal (i.e., equality holds above),
we must have $S X' = 0$; we use this fact later.

To show that $Y$ is a globally optimal solution to \eqref{eq:rankp_laplace_sbd},
it clearly suffices to show that $Y Y^\top$ is an optimal solution to \eqref{eq:sdp_laplace_sbd}.
Thus we want to find a dual certificate proving optimality of $Y Y^\top$.
A natural candidate is the matrix $S(Y) = \Lhat - \sbd(\Lhat Y Y^\top)$ from \eqref{eq:SY}.
First-order criticality implies that $S(Y) Y Y^\top = 0$.
Thus it remains to show that $S(Y) \succeq 0$.
We do this by leveraging the following fact (see, e.g., \cite[Prop.~3.1]{Boumal2019}):
\begin{lemma}
	\label{lem:rankdef}
	If $Y$ is a second-order critical point of \eqref{eq:rankp_laplace_sbd} that is rank deficient (i.e., $\rank(Y) < p$),
	then $S(Y) \succeq 0$, and consequently $Y Y^\top$ is an optimal solution to the SDP \eqref{eq:sdp_laplace_sbd},
	and $Y$ is a globally optimal solution to \eqref{eq:rankp_laplace_sbd}.
\end{lemma}
This can be proved with an argument from~\cite{Journee2010}:
if $u$ is any unit-norm vector in the null space of $Y$, then, for any $z \in \R^{rn}$,
$\Ydt(z) \coloneqq z u^\top$ is in $T_Y$; hence, second-order criticality implies
\[
	\ip{S(Y) z}{z} = \ip{S(Y) \Ydt(z)}{\Ydt(z)} \geq 0.
\]
Thus the key step in showing that $Y$ is globally optimal is to show that it is \emph{rank deficient}.

\section{Proofs (real case)}
\label{sec:real_proofs}

We (still) denote by $\opnorm{A}$, $\normF{A}$, and $\nucnorm{A}$ the operator, Frobenius, and nuclear norms of a matrix $A$.
Our proof strategy is as follows.

First, we show that each second-order critical point $Y$ of \eqref{eq:rankp_orig} satisfies certain error bounds (\Cref{thm:noiseless,thm:noisy_bound}).
To do this, we can use inequalities of the form $\ip{S(Y) \Ydt}{\Ydt} \geq 0$ for all $\Ydt \in T_Y$---these express the fact that the Riemannian Hessian at $Y$ is positive semidefinite, with $S(Y)$ as in~\eqref{eq:SY}.
We did not find a way to select a single $\Ydt$ to prove our claims.
Instead, the key enabling proof idea is this:
we design a suitable probability distribution over the tangent space $T_Y$,
and we exploit the fact that $\ip{S(Y) \Ydt}{\Ydt}$ is still nonnegative \emph{in expectation} over the random choice of $\Ydt$.
This yields a valid inequality which may not be expressible using a single $\Ydt$.

Second, to prove our landscape results (\Cref{thm:noisy_rankbd,thm:noisy_landscape}),
the critical step is to show that second-order critical points have low rank (after which we apply \Cref{lem:rankdef}).
To do this, we show (using the error bounds previously derived) that the matrix $S(Y)$ has \emph{high} rank and then apply the first-order criticality condition $S(Y) Y = 0$.

To simplify notation in our proofs, we assume that $Z_i = I_r$ for all $i = 1, \dots, n$.
This is without loss of generality because it can be arranged with a smooth and smoothly invertible change of variable.
More explicitly, if $Z$ is arbitrary (with $Z_i \in \Or$), consider the change of variable $\pi \colon \scrM \to \scrM$ defined by $\pi(Y)_i = Z_i Y_i$.
Since $\pi$ is a Riemannian isometry, the landscapes of $F(Y) = \ip{\Lhat}{YY^\top}$ \eqref{eq:rankp_laplace_sbd} and of $\Fbr = F\circ\pi$ are the same \cite[Prop.~9.6]{Boumal2023}, in the sense that
$Y$ is first-order critical / second-order critical / locally optimal / globally optimal for $\Fbr$ if and only if $\pi(Y)$ is so for $F$.
Using the definition of $\Lhat$, the expression for $\Fbr$ simplifies to $\Fbr(Y) = \ip{L \otimes I_r - \Deltabr}{YY^\top}$, where $\Deltabr_{ij} = Z_i^\top\Delta_{ij}Z_j$.
This is exactly $F$ in the event that $Z_1 = \cdots = Z_n = I_r$, including a change of variable on the noise matrix $\Delta \mapsto \Deltabr$ which has no effect on its eigenvalues (hence on any of the claims we make about noise later on).

Accordingly, we proceed with the following simplified notation:
\begin{align*}
    Z_i = I_r\ \forall i, &&
    D = I_{rn}, &&
    L_Z = L \otimes I_r, &&
    \Lhat = L_Z - \Delta.
\end{align*}
The next step is to design a distribution of random tangent vectors.

\subsection{A distribution of random tangent vectors}
\label{sec:rand_tangent}
We analyze critical points of \eqref{eq:rankp_laplace_sbd}.
Recall from \Cref{sec:math_background} that, setting $S(Y) = \Lhat - \sbd(\Lhat Y Y^\top)$,
a feasible point $Y$ is first- and second-order critical if $S(Y) Y = 0$ and, for every $\Ydt \in T_Y$, $\ip{S(Y)}{\Ydt \Ydt^\top} \geq 0$.
From here on, assume $Y$ is such a point.

Recall that $\Ydt$ is in $T_Y$ if and only if $\Ydt_i^{} Y_i^\top + Y_i^{} \Ydt_i^\top = 0$ for $i = 1, \dots, n$.
In other words, each $Y_i^{} \Ydt_i^\top$ must be skew-symmetric.
This is true if and only we can write
\[
	\Ydt_i = \Gamma_i (I - Y_i^\top Y_i^{}) + S_i Y_i
\]
for some $r \times p$ matrix $\Gamma_i$ and some \emph{skew-symmetric} $r \times r$ matrix $S_i$.
The first term is the row-wise orthogonal projection of $\Gamma_i$ onto $\nul(Y_i)$.
If we choose $S_i = \Gamma_i^{} Y_i^\top - Y_i^{} \Gamma_i^\top$,
we obtain
\[
	\Ydt_i = \Gamma_i - Y_i^{} \Gamma_i^\top Y_i^{} = Y_i (Y_i^\top \Gamma_i^{} - \Gamma_i^\top Y_i^{}).
\]
In fact, this last formulation covers all possible $\Ydt_i$,
because any $r \times r$ skew-symmetric matrix $S_i$ can be written in the given form,
and the components of $\Gamma_i$ in $\range(Y_i^\top)$ and $\nul(Y_i)$ can be chosen independently.

We choose a common $\Gamma_i = \Gamma$, where $\Gamma$ is a random $r \times p$ matrix whose entries are i.i.d.\ standard normal random variables.
This results in a random $\Ydt \in T_Y$ whose components $\Ydt_1, \ldots, \Ydt_n$ are related.
Because the second-order criticality inequality holds for \emph{all} $\Ydt \in T_Y$, we can take an expectation to obtain
\[
	\ip{S(Y)}{\E \Ydt \Ydt^\top} = \E \ip{S(Y)}{\Ydt \Ydt^\top} \geq 0.
\]

Simple calculations yield
\begin{align}
    \E \Ydt_i^{} \Ydt_j^\top
        & = \E ( \Gamma - Y_i \Gamma^\top Y_i )(\Gamma - Y_j \Gamma^\top Y_j)^\top \nonumber\\
        & = \E \parens*{ \Gamma \Gamma^\top - Y_i \Gamma^\top Y_i \Gamma^\top - \Gamma Y_j^\top \Gamma Y_j^\top + Y_i^{} \Gamma^\top Y_i^{} Y_j^\top \Gamma Y_j^\top } \nonumber\\
        & = (p - 2) I_r + \tr(Y_i^{} Y_j^\top) Y_i^{} Y_j^\top.
    \label{eq:tangent_cov_real}
\end{align}
We have used the facts that (1) for any $r \times p$ matrix $U$, $\E \Gamma^\top U \Gamma^\top = U^\top$,
and (2) for any $r \times r$ matrix $B$, $\E \Gamma^\top B \Gamma = \tr(B) I_p$.
When $i = j$, this simplifies to $\E \Ydt_i^{} \Ydt_i^\top = (p + r - 2) I_r$. This implies $\sbd(\E \Ydt \Ydt^\top) = (p + r - 2) I_{rn}$.

Using a random $\Ydt \in T_Y$ has appeared before in \cite{Mei2017,Ling2022b}.
Those papers chose (in our notation) $\Ydt_i = \Gamma(I - Y_i^\top Y_i)$, which is only the portion of $\Gamma$ in $\nul(Y_i)$.
Another similar approach to ours appears in \cite{Markdahl2018a,Markdahl2020},
where, rather than considering a random choice of $\Ydt$,
the authors analyze the quadratic form $\Gamma \mapsto \ip{S(Y)}{\Ydt(\Gamma) \Ydt(\Gamma)^\top}$ with $\Ydt_i(\Gamma) = \scrP_{T_{Y_i}}(\Gamma)$.
Their results do not use randomness and scale differently the portions of each $\Ydt_i$ that are in $\nul(Y_i)$ and $\range(Y_i^\top)$ (this is related to the choice of metric on the Stiefel manifold).
\subsection{Noiseless case}
\label{sec:proof_noiseless}
For \Cref{thm:noiseless}, we have $\Delta = 0$, so $\Lhat = L_Z = L \otimes I_r$, and $S(Y) = L_Z - \sbd(L_Z Y Y^\top)$.

We first calculate
\begin{align*}
	\ip{\sbd(\Lhat Y Y^\top)}{\E \Ydt \Ydt^\top}
	&= \ip{L_Z Y Y^\top}{\sbd(\E \Ydt \Ydt^\top)} \\
	&= (p + r - 2) \tr(L_Z Y Y^\top) \\
	&= (p + r - 2) \sum_{i,j=1}^n L_{ij} \tr(Y_i^{} Y_j^\top).
\end{align*}

Next, by \eqref{eq:tangent_cov_real}, we derive
\begin{align*}
	\ip{L_Z}{\E \Ydt \Ydt^\top}
	&= \sum_{i,j=1}^n L_{ij} \ip{I_r}{ \E \Ydt_i^{} \Ydt_j^\top } \\
	&= \sum_{i,j=1}^n L_{ij} ((p-2)r + \tr^2(Y_i^{} Y_j^\top) ) \\
	&= \sum_{i,j=1}^n L_{ij} \tr^2(Y_i^{} Y_j^\top),
\end{align*}
where the last equality follows from the fact that $\sum_i L_{ij} = 0$ for all $j$.

To proceed, note that
\[
	Y_i^{} Y_j^\top + Y_j^{} Y_i^\top
	= 2 I_r - (Y_i - Y_j)(Y_i - Y_j)^\top,
\]
so that
\begin{equation}
	\label{eq:trace_id_real}
	\tr(Y_i^{} Y_j^\top) = \frac{1}{2} \tr(Y_i^{} Y_j^\top + Y_j^{} Y_i^\top) = r - \frac{1}{2} \normF{Y_i - Y_j}^2.
\end{equation}
We then find (again using the fact that all rows and columns of $L$ sum to zero)
\begin{align*}
	\ip{L_Z}{\E \Ydt \Ydt^\top}
	&= \sum_{i,j=1}^n L_{ij} \parens*{ r - \frac{1}{2}\normF{Y_i - Y_j}^2 } ^2 \\
	&= -r \sum_{i,j=1}^n L_{ij} \normF{Y_i - Y_j}^2 + \frac{1}{4} \sum_{i,j=1}^n L_{ij}  \normF{Y_i - Y_j}^4 \\
	&= 2 r \sum_{i,j=1}^n L_{ij} \tr(Y_i^{} Y_j^\top) - \frac{1}{4} \sum_{i,j=1}^n A_{ij} \normF{Y_i - Y_j}^4.
\end{align*}

All combined, the condition $\ip{S(Y)}{ \E \Ydt \Ydt^\top} \geq 0$ then implies
\begin{equation}
	\label{eq:final_ineq_noiseless}
	(p - r - 2) \sum_{i,j=1}^n L_{ij} \tr(Y_i^{} Y_j^\top) + \frac{1}{4} \sum_{i,j=1}^n A_{ij} \normF{Y_i - Y_j}^4 \leq 0.
\end{equation}
The second term on the left-hand side is nonnegative. If $p \geq r + 2$, the first term is also nonnegative (the sum equals $\ip{L_Z}{Y Y^\top}$).
Therefore, both terms on the left-hand side of \eqref{eq:final_ineq_noiseless} must be $0$, so $Y_i = Y_j$ for all $(i,j) \in E$.
If $G$ is connected, this implies that the $Y_i$'s are identical, and therefore
\[
	\ip{Z Z^\top}{Y Y^\top} = \normF{Z^\top Y}^2 = \normF{n Y_1}^2 = n^2 r.
\]
This finishes the proof of \Cref{thm:noiseless}.

When $p = r + 2$, the robustness to (small) perturbations comes from the quartic terms in \eqref{eq:final_ineq_noiseless}.
This robustness is much weaker than what we prove in \Cref{sec:proof_noisy} for the noisy results (where we use the \emph{quadratic} terms that arise when $p > r + 2$).

\subsection{Noisy case}
\label{sec:proof_noisy}
We now consider the case where the noise matrix $\Delta$ is nonzero.
We prove \Cref{thm:noisy_bound} first (high correlation) and then use it to prove \Cref{thm:noisy_rankbd,thm:noisy_landscape} (global optimality).

To aid our analysis, we write $Y = Z R + W$, with $R = \frac{1}{n} Z^\top Y \in \R^{r \times p}$ and $W$ orthogonal to $Z$ (i.e., $Z^\top W = \sum_i W_i = 0$). This gives $Y_i = R + W_i$.
Note that
\begin{equation}
	\label{eq:error_ids}
	\normF{Z^\top Y}^2 = \normF{nR}^2 = n \normF{Z R}^2 = n^2 r - n \normF{W}^2.
\end{equation}
Thus, we set out to use second-order criticality of $Y$ to infer that $\normF{W}$ is small.

\subsubsection{Error bound for all second-order critical points}
The second-order criticality condition is now, for all $\Ydt \in T_Y$,
\begin{equation}
	\label{eq:soc_noisy}
	\ip{L_Z - \Delta}{\Ydt \Ydt^\top} - \ip{\sbd((L_Z - \Delta) Y Y^\top)}{\Ydt \Ydt^\top} \geq 0.
\end{equation}
We use the same random $\Ydt$ as in the previous section.
We have already analyzed the terms involving $L_Z$ in \Cref{sec:proof_noiseless}, so it remains to analyze the terms involving $\Delta$.

Note that because we have assumed that every $Z_i = I_r$, we have
\begin{equation}
	\label{eq:Delta_Z_ip}
	\sum_{i,j=1}^n \ip{\Delta_{ij}}{I_r} = \sum_{i,j=1}^n \ip{\Delta_{ij}}{Z_i^{} Z_j^\top} = \ip{\Delta}{Z Z^\top}.
\end{equation}
From \eqref{eq:tangent_cov_real}, we calculate
\begin{align*}
	\ip{\Delta}{\E \Ydt \Ydt^\top}
	&= \sum_{i,j = 1}^n \ip{\Delta_{ij}}{(p-2) I_r + \tr(Y_i^{} Y_j^\top) Y_i^{} Y_j^\top} \\
	&= (p-2)\ip{\Delta}{Z Z^\top} + \sum_{i,j=1}^n \parens*{ r - \frac{1}{2} \normF{Y_i - Y_j}^2 } \ip{\Delta_{ij}}{Y_i^{} Y_j^\top} \\
	&= (p-2)\ip{\Delta}{Z Z^\top} + r \ip{\Delta}{Y Y^\top} - \frac{1}{2} \sum_{i,j = 1}^n \normF{Y_i - Y_j}^2 \ip{\Delta_{ij}}{Y_i^{} Y_j^\top}.
\end{align*}
The second equality uses \eqref{eq:trace_id_real} and \eqref{eq:Delta_Z_ip}.
Next, $\sbd(\E \Ydt \Ydt^\top) = (p + r - 2) I_{rn}$ implies
\[
	\ip{\sbd(\Delta Y Y^\top)}{\E \Ydt \Ydt^\top} = (p + r - 2) \ip{\Delta}{Y Y^\top}.
\]
The difference is
\begin{align}
	&\negqquad \ip{\sbd(\Delta Y Y^\top)}{\E \Ydt \Ydt^\top} - \ip{\Delta}{\E \Ydt \Ydt^\top} \nonumber\\
	&= (p - 2) \ip{\Delta}{Y Y^\top - Z Z^\top} + \sum_{i,j = 1}^n \ip*{\Delta_{ij}}{\frac{1}{2} \normF{Y_i - Y_j}^2 Y_i^{} Y_j^\top} \nonumber\\
	&= (p - 2) \ip{\Delta}{Y Y^\top - Z Z^\top} + \ip{\Delta}{(Q \otimes \onevec_{r}^{} \onevec_r^\top) \circ (Y Y^\top)},
	\label{eq:gohereforholder}
\end{align}
where $\circ$ is entrywise product and $Q_{ij} = \frac{1}{2} \normF{Y_i - Y_j}^2$.

We aim to bound (in nuclear norm) the matrices that are in an inner product with $\Delta$.
Note that the $r \times r$ blocks of $Y Y^\top - Z Z^\top$ are
\[
	(Y Y^\top - Z Z^\top)_{ij} = Y_i Y_j^\top - I_r = Y_i(Y_j - Y_i)^\top = Y_i(W_j - W_i)^\top,
\]
so we can write $Y Y^\top - Z Z^\top = Y W^\top - H$,
where $H_{ij} = Y_i^{} W_i^\top$.
By the matrix Hölder inequality for Schatten $p$-norms, $\nucnorm{Y W^\top} \leq \normF{Y} \normF{W} = \sqrt{rn} \normF{W}$.
Furthermore,
\[
H = \begin{bmatrix*}
	Y_1^{} W_1^\top \\ \vdots \\ Y_n^{} W_n^\top
\end{bmatrix*} Z^\top.
\]
Because each $Y_i$ has operator norm $1$, the left factor in the above expression has Frobenius norm at most $\normF{W}$.
Also, $\normF{Z} = \sqrt{rn}$.
Thus $\nucnorm{H} \leq \sqrt{rn} \normF{W}$. We conclude that $\nucnorm{Y Y^\top - Z Z^\top} \leq 2 \sqrt{nr} \normF{W}$.

We also need a bound on
\begin{align*}
	\nucnorm{(Q \otimes \onevec_{r}^{} \onevec_r^\top) \circ (Y Y^\top)}
	\overset{\mathrm{(a)}}{\leq} \nucnorm{(Q \otimes \onevec_{r}^{} \onevec_r^\top)}
	\overset{\mathrm{(b)}}{=} r \nucnorm{Q}. 
\end{align*}
Inequality (a) follows from a basic inequality on singular values of Hadamard products \cite[Thm.~5.6.2]{Horn1991}, using the fact that each row of $Y$ has norm $1$.
Equality (b) follows from the eigenvalue characterization of Kronecker products.
To bound $\nucnorm{Q}$, note that $Q_{ij}$ is simply the trace of the $(i,j)$th block of $Z Z^\top - Y Y^\top$.
Thus $Q$ is a partial trace of $Z Z^\top - Y Y^\top$,
and therefore (see \cite{Rastegin2012}) $\nucnorm{Q} \leq \nucnorm{Z Z^\top - Y Y^\top} \leq 2 \sqrt{nr} \normF{W}$.

Putting all the nuclear norm bounds together, we obtain, by Hölder's inequality for matrix inner products applied on~\eqref{eq:gohereforholder} (von Neumann's trace inequality),
\[
	\ip{\sbd(\Delta Y Y^\top)}{\E \Ydt \Ydt^\top} - \ip{\Delta}{\E \Ydt \Ydt^\top}
	\leq 2(p + r - 2) \opnorm{\Delta} \sqrt{rn}\normF{W}.
\]

Combining this with \eqref{eq:soc_noisy} and the calculations of \Cref{sec:proof_noiseless},
we obtain
\begin{equation}
	\label{eq:LYY_bound}
	(p - r - 2)\sum_{i,j} L_{ij} \ip{Y_i}{Y_j} + \frac{1}{4} \sum_{i,j=1}^n A_{ij} \normF{W_i - W_j}^4
	\leq 2(p + r - 2) \opnorm{\Delta} \sqrt{rn}\normF{W}.
\end{equation}
Dropping the nonnegative quartic terms and using the fact that
\[
	\sum_{i,j} L_{ij} \ip{Y_i}{Y_j} = \sum_{i,j} L_{ij} \ip{W_i}{W_j} \geq \lambda_2 \normF{W}^2
\]
(where $\lambda_2 = \lambda_2(L)$ is the Fiedler value of $G$), we get
\begin{align*}
	(p - r - 2) &\lambda_2 \normF{W}^2 \leq 2(p + r - 2) \opnorm{\Delta} \sqrt{rn}\normF{W}.
\end{align*}
Therefore,
\begin{align}
	\normF{W}^2 \leq \parens*{ \frac{2(p + r - 2) \opnorm{\Delta}}{(p - r - 2)\lambda_2}}^2 rn
	= C_p^2 rn \frac{\opnorm{\Delta}^2}{\lambda_2^2}.
    \label{eq:boundonWF}
\end{align}
Combining this with the identities \eqref{eq:error_ids} finishes the proof of \Cref{thm:noisy_bound}.

\subsubsection{Bound on solution rank}
\label{sec:proof_noisy_rankbound}
We next prove \Cref{thm:noisy_rankbd}.
To show that a second-order critical point $Y$ has low rank,
recall that the first-order criticality condition states $S(Y) Y = 0$.
Therefore, it suffices to show that $S(Y)$ has high rank.
Recall
\[
	S(Y) = \Lhat - \sbd(\Lhat Y Y^\top) = L_Z - \Delta - \sbd(L_Z Y Y^\top) + \sbd(\Delta Y Y^\top).
\]
Because $G$ is connected, $\lambda_2$ is positive.
Thus $L_Z$ has an $r$-dimensional null space,
and its remaining $(n-1)r$ eigenvalues are at least $\lambda_2$.
Assume $\opnorm{\Delta} < \lambda_2$ (as otherwise the theorem statement is true but vacuous).
Given a matrix $A$, let $\sigma_k(A)$ denote its $k$th singular value, in decreasing order.
To bound the rank of $Y$, we count how many singular values of $S(Y)$ can possibly be equal to zero.
Explicitly, for any $c \in (0, 1)$, it holds (see the explanation below):
\begin{align}
	\rank(Y)
	&\leq \dim \nul(S(Y)) \nonumber\\
	&\leq r + \abs{ \{ \ell : \sigma_\ell(\sbd(L_Z Y Y^\top) -  \sbd(\Delta Y Y^\top)) \geq \lambda_2 - \opnorm{\Delta} \} } \nonumber\\
	&\leq r + \abs{ \{ \ell: \sigma_\ell(\sbd(L_Z Y Y^\top)) \geq  c (\lambda_2 - \opnorm{\Delta}) \} } \nonumber\\
	&\qquad + \abs{ \{ \ell: \sigma_\ell(\sbd(\Delta Y Y^\top)) \geq  (1-c) (\lambda_2 - \opnorm{\Delta}) \} } \nonumber\\
	&\leq r + \frac{\nucnorm{\sbd(L_Z Y Y^\top)}}{c(\lambda_2 - \opnorm{\Delta})} + \frac{\normF{\sbd(\Delta Y Y^\top)}^2}{(1 - c)^2 (\lambda_2 - \opnorm{\Delta})^2}.
    \label{eq:rankY_bd_first}
\end{align}
The third inequality is an application of the following fact: if matrices $A$ and $B$ (of the same size) have, respectively, at most $k_A$ singular values larger than $M_A$ and $k_B$ singular values larger than $M_B$, then $A + B$ has at most $k_A + k_B$ singular values larger than $M_A + M_B$.
This follows from \cite[Thm.~3.3.16]{Horn1991}, which implies $\sigma_{k_A + k_B + 1}(A+B) \leq \sigma_{k_A+1}(A) + \sigma_{k_B+1}(B) < M_A + M_B$.

To bound this last quantity, first note that, because $\opnorm{Y_i} = 1$ for each $i$,
\begin{equation}
	\label{eq:DYY_fronorm_bd}
	\normF{\sbd(\Delta Y Y^\top)}^2 \leq \sum_i \normF{(\Delta Y)_i^{} Y_i^\top}^2 \leq \normF{\Delta Y}^2 \leq \opnorm{\Delta}^2 \normF{Y}^2 = rn \opnorm{\Delta}^2.
\end{equation}

Next, using $L_Z = L \otimes I_r$, note that
\begin{align*}
	(\sbd(L_Z Y Y^\top))_{ii}
	&= \frac{1}{2} \sum_{j=1}^n L_{ij} (Y_i^{} Y_j^\top + Y_j^{} Y_i^\top) \\
	&= \sum_{j=1}^n L_{ij} \parens*{I_r - \frac{1}{2} (Y_i - Y_j)(Y_i - Y_j)^\top} \\
	&= \frac{1}{2} \sum_{j=1}^n A_{ij} (Y_i - Y_j)(Y_i - Y_j)^\top \\
	&\succeq 0.
\end{align*}
Therefore, $\sbd(L_Z Y Y^\top)$ is positive semidefinite and it follows that its nuclear norm is bounded as:
\begin{equation}
	\label{eq:LYY_nucnorm_bd}
\begin{aligned}
	\nucnorm{\sbd(L_Z Y Y^\top)}
	&= \tr(\sbd(L_Z Y Y^\top)) \\
	&= \sum_{i,j=1}^n L_{ij} \tr(Y_i^{} Y_j^\top) \\
	&\overset{\text{(i)}}{\leq} C_p \opnorm{\Delta} \sqrt{rn} \normF{W} \\
	&\overset{\text{(ii)}}{\leq} C_p^2 rn \frac{\opnorm{\Delta}^2}{\lambda_2}.
\end{aligned}
\end{equation}
Inequality (i) comes from \eqref{eq:LYY_bound} in the proof of \Cref{thm:noisy_bound}.
Inequality (ii) follows from the bound on $\normF{W}$ provided by \eqref{eq:boundonWF} in that same proof.

Plugging \eqref{eq:DYY_fronorm_bd} and \eqref{eq:LYY_nucnorm_bd} into \eqref{eq:rankY_bd_first},
we obtain
\begin{align*}
	\rank(Y)
	&\leq r + \parens*{ \frac{C_p^2}{c \lambda_2 (\lambda_2 - \opnorm{\Delta})} + \frac{1}{(1 - c)^2 (\lambda_2 - \opnorm{\Delta})^2} } \opnorm{\Delta}^2 rn.
\end{align*}
Assume $\opnorm{\Delta} \leq \lambda_2 / 4$.
Choosing $c = 1/2$ and using the fact that $C_p \geq 2$, 
we see that
\begin{align*}
	\rank(Y) 
	&\leq r + \parens*{ \frac{8}{3} \frac{ C_p^2}{\lambda_2^2} + \frac{64}{9} \frac{1}{\lambda_2^2} } \opnorm{\Delta}^2 rn
	\leq r + 5 \parens*{ \frac{C_p \opnorm{\Delta}}{\lambda_2}}^2 rn.
\end{align*}
That last bound still holds if $\opnorm{\Delta} > \lambda_2 / 4$ since $\rank(Y) \leq rn$.
This completes the rank-bound portion of \Cref{thm:noisy_rankbd}.
If $p$ is larger than this rank bound, then $Y$ is rank-deficient,
and we apply \Cref{lem:rankdef} to obtain the rest of the result.

A limiting factor in how small we must make $\opnorm{\Delta}$ in the above proof is the need to bound the singular values of $\sbd(L_Z Y Y^\top)$.
If $G$ is the complete graph, $L_Z = I_{rn} - Z Z^\top$,
so
\[
	S(Y) = \sbd((Z Z^\top + \Delta) Y Y^\top) - Z Z^\top - \Delta.
\]
Here it suffices to \emph{lower bound} the singular values of the first term,
and there is no need to bound $\sbd(L_Z Y Y^\top)$.
This simplification allows existing work (e.g., \cite{Ling2022b}) to obtain better results in the complete-graph case (with additional assumptions on $\Delta$).
It is not clear how to analyze a general graph $G$
in such a way that we recover existing results when $G$ is complete.
\subsubsection{Solution uniqueness}
To prove \Cref{thm:noisy_landscape},
note that, under the assumptions of \Cref{thm:noisy_landscape}, \Cref{thm:noisy_rankbd} and its proof imply that $\rank(S(Y)) \geq rn - r$ and, by first-order criticality, $\rank(Y) \leq r$.
These rank inequalities are, in fact, equalities, because the constraint $Y_1^{} Y_1^\top = I_r$ implies $\rank(Y) \geq r$.
Thus $\rank(S(Y)) = rn - r$, and $\rank(Y) = r$.

Since $p > r$, \Cref{lem:rankdef} again implies that $Y Y^\top$ solves the SDP \eqref{eq:sdp_laplace_sbd},
and $S(Y) \succeq 0$ is its dual certificate.
We set out to prove that $Y Y^\top$ is the \emph{unique} solution.
First, note that because $Y$ has rank $r$ we can write $Y = \Zhat U$, with $\Zhat \in \Ogrp(r)^n$ and $U \in \R^{r \times p}$ such that $U U^\top = I_r$.
Then $Y Y^\top = \Zhat \Zhat^\top$.
Furthermore, as discussed in \Cref{sec:dualcerts}, \emph{any} optimal solution $X$ of the SDP must satisfy $S(Y)X = 0$.
Because the columns of $\Zhat$ span the kernel of $S(Y)$,
we can write $X = (\Zhat B)(\Zhat B)^\top = \Zhat (B B^\top) \Zhat^\top$ for some $r \times r$ matrix $B$.
The block-diagonal constraint implies that
\[
	0 = (Y Y^\top - X)_{11} = (\Zhat \Zhat^\top - \Zhat (B B^\top)\Zhat^\top)_{11}
	= Z_1 (I_r - B B^\top) Z_1^\top.
\]
Because $Z_1 \in \Ogrp(r)$, we must have $B B^\top = I_r$, and therefore $X = Y Y^\top$.
Thus $Y Y^\top = \Zhat \Zhat^\top$ is the \emph{unique} SDP solution.

The fact that $Y Y^\top = \Zhat \Zhat^\top$ is the optimal solution to \eqref{eq:sdp_laplace_sbd} (and therefore of \eqref{eq:sdp_orig}) implies that $\Zhat$ and $Y$ are, respectively, global optima of \eqref{eq:ncvx_orig} and \eqref{eq:rankp_orig}.
For uniqueness, note that for any other global optima $Z'$ of \eqref{eq:ncvx_orig} and $Y'$ of \eqref{eq:rankp_orig},
$Z' Z'^\top$ and $Y' Y'^\top$ are feasible points of \eqref{eq:sdp_orig} with the same objective function value as $\Zhat \Zhat^\top = Y Y^\top$. Therefore, by the uniqueness of the SDP solution, we have $Z' Z'^\top = Y' Y'^\top = Y Y^\top = \Zhat \Zhat^\top$,
implying that $Z' = \Zhat$ and $Y' = Y$ up to global orthogonal transformations.
This completes the proof of \Cref{thm:noisy_landscape}.

\section{Extension of proofs to the complex case}
\label{sec:complexcase}
We extend the previous section's argument to the complex case in a direct way,
only considering noiseless measurements.
This yields the apparently new yet possibly suboptimal result stated in \Cref{thm:complex}.
We simply substitute complex quantities for real ones in the previous arguments,
replacing transposes ($A^\top$) by their Hermitian counterparts ($A^*$).
To avoid ambiguity, we will define the matrix inner product as $\ip{A}{B} = \real(\tr(A B^*))$,
writing out the trace explicitly when we need to consider its imaginary part.

We still assume, without loss of generality, that $Z_1 = \cdots = Z_n = I_r$, so that $\Lhat = L_Z = L \otimes I_r$ in the noiseless case.
Once again, we use the fact that the Riemannian Hessian is PSD at a second-order critical point.
Similarly to \eqref{eq:SY}, we now take $S(Y) = \Lhat - \sbd(\Lhat Y Y^*)$, where $\sbd$ now keeps the \emph{Hermitian} part of the $r\times r$ diagonal blocks (again making the off-diagonal blocks zero).
Second-order criticality means that $S(Y)Y = 0$ and $\ip{S(Y)}{\Ydt \Ydt^*} \geq 0$ for all $\Ydt$ in the tangent space at $Y \in \St(r, p, \C)^n$.
Details appear in~\cite{Pumir2018}.

For the time being (we shall modify this later in this section), we directly translate the construction in \Cref{sec:rand_tangent} and choose tangent vectors $\Ydt_i = \Gamma - Y_i \Gamma^* Y_i$ to the points $Y_i \in \St(p, r, \C)$,
where now $\Gamma$ is an $r \times p$ matrix of i.i.d.\ \emph{complex} standard normal random variables.
This is indeed tangent as $\Ydt_i^{} Y_i^* + Y_i^{} \Ydt_i^* = 0$.
By a similar calculation to before,
\begin{align*}
	\E \Ydt_i^{} \Ydt_j^* &= \E ( \Gamma \Gamma^* - Y_i \Gamma^* Y_i \Gamma^* - \Gamma Y_j^* \Gamma Y_j^* + Y_i^{} \Gamma^* Y_i^{} Y_j^* \Gamma Y_j^* ) \\
	&= p I_r + \tr(Y_i^{} Y_j^*) Y_i^{} Y_j^*.
\end{align*}
The first difference from the real case is that the terms with two factors of $\Gamma$ or $\Gamma^*$ have zero expectation, because each entry in these terms is a polynomial in standard complex normal random variables and thus is radially symmetric.\footnote{If $z$ is a standard normal random variable (scalar), then in the real case $\E z^2 = 1$ but in the complex case $\E z^2 = 0$.}
Note that $\E \Ydt_i^{} \Ydt_i^* = (p+r) I_r$.
First, we can compute, similarly to before,
\begin{align}
	\ip{\sbd(L_Z Y Y^*)}{\E \Ydt \Ydt^*}
	&= (p + r) \sum_{i,j=1}^n L_{ij} \ip{Y_i}{Y_j}.
    \label{eq:sbdcomplexineq}
\end{align}

Next, we have (see \Cref{sec:proof_noiseless})
\begin{align}
   	\ip{L_Z}{\E \Ydt \Ydt^*} = \sum_{i,j=1}^n \real(\tr^2(Y_i^{} Y_j^*) ).
\end{align}
We now come to the second significant difference from the real case:
the calculation \eqref{eq:trace_id_real} fails here,
because $\tr(Y_i^{} Y_j^*)$ is not necessarily real.\footnote{If we ignored this issue, we could immediately ``prove'' that $p = r$ suffices to obtain a benign landscape. This is clearly false, as the case $r = 1$ (angular synchronization) demonstrates.}
By considering the real and imaginary parts of $\tr (Y_i^{} Y_j^*)$, we obtain
\begin{align}
    \ip{L_Z}{\E \Ydt \Ydt^*} = \frac{1}{4} \sum_{i,j=1}^n L_{ij} \brackets*{ \tr^2(Y_i^{} Y_j^* + Y_j^{} Y_i^*) - \abs{\tr(Y_i^{} Y_j^* - Y_j^{} Y_i^*)}^2 }.
\end{align}
The first term can be handled akin to \eqref{eq:trace_id_real} in the real case.
The optimal way to handle the second term is unclear.
One way is to use $|\tr(A)|^2 \leq \nucnorm{A}^2 \leq r \normF{A}^2$ (twice) and $\normF{A+A^*}^2 + \normF{A-A^*}^2 = 4\normF{A}^2$ to show that
\begin{equation}
	\label{eq:imag_bound}
	\begin{aligned}
		\abs{\tr(Y_i^{} Y_j^* - Y_j^{} Y_i^*)}^2
		&\leq r \normF{Y_i^{} Y_j^* - Y_j^{} Y_i^*}^2 \\
		&= 4 r \normF{Y_i^{} Y_j^*}^2 - r \normF{Y_i^{} Y_j^* + Y_j^{} Y_i^*}^2 \\
		&\leq 4 r \normF{Y_i^{} Y_j^*}^2 - \tr^2( Y_i^{} Y_j^* + Y_j^{} Y_i^* ).
	\end{aligned}
\end{equation}
Both sides of the above inequality are zero when $i = j$.
For $i \neq j$, recall that $L_{ij} \leq 0$.
Therefore,
\begin{align}
    \ip{L_Z}{\E \Ydt \Ydt^*} \leq \frac{1}{4} \sum_{i,j=1}^n L_{ij} \brackets*{ 2\tr^2(Y_i^{} Y_j^* + Y_j^{} Y_i^*) - 4 r \normF{Y_i^{} Y_j^*}^2 }.
\end{align}
Additionally, $\sum_{ij} L_{ij} \normF{Y_i^{} Y_j^*}^2 = \sum_{ij} L_{ij} \ip{Y_i^* Y_i^{}}{Y_j^* Y_j^{}} \geq 0$, hence we can remove the last term while preserving the inequality.
We now combine this and~\eqref{eq:sbdcomplexineq} into the inequality $\ip{S(Y)}{\E \Ydt \Ydt^*} \geq 0$ to obtain
\begin{align*}
	(p + r) \sum_{i,j=1}^n L_{ij} \ip{Y_i}{Y_j}
	&\leq \frac{1}{2} \sum_{i,j=1}^n L_{ij} \tr^2( Y_i^{} Y_j^* + Y_j^{} Y_i^* ) \\
	&= 4r \sum_{i,j=1}^n L_{ij} \ip{Y_i}{Y_j} - \frac{1}{2} \sum_{i,j=1}^n A_{ij} \normF{Y_i - Y_j}^4.
\end{align*}
(The equality follows just as in the real case, e.g., by \eqref{eq:trace_id_real}.)

This latter inequality yields perfect recovery provided $p \geq 3r$.
This condition seems much worse than what we obtained in the real case.

We can improve the result by making a slightly different choice of the $\Ydt_i$'s.
Note that we can write our choice as
\[
	\Ydt_i = \Gamma - Y_i \Gamma^* Y_i
	= \Gamma(I_p - Y_i^* Y_i^{}) + (\Gamma Y_i^* - Y_i^{} \Gamma^*) Y_i.
\]
The first term has rows orthogonal to $Y_i$, while the second has a skew-symmetric matrix left-multiplying $Y_i$ (rather than right-multiplying as before): both terms are therefore tangent vectors.
We can rescale them arbitrarily, as
\[
	\Ydt_i = a \Gamma(I_p - Y_i^* Y_i^{}) + b (\Gamma Y_i^* - Y_i^{} \Gamma^*) Y_i
\]
for any numbers $a \in \C$ and $b \in \R$.
For $a, b > 0$ (which turns out to be the only sensible choice), this is related to the choice of metric on the Stiefel manifold.
We choose $a = 2$ and $b = 1$, which makes $\Ydt$ the (Euclidean) orthogonal projection of $2 Z \Gamma$ onto $T_Y$.
Intrinsically, this is quite similar to a complex adaptation of the proof in \cite{Markdahl2020}, though that argument is not phrased in terms of randomness.
Now, for $\Gamma$ chosen randomly as before,
\begin{align*}
	\E \Ydt_i^{} \Ydt_j^*
	&= 4 \E \Gamma(I_p - Y_i^* Y_i^{}) (I_p - Y_j^* Y_j^{}) \Gamma^* + \E (\Gamma Y_i^* - Y_i^{} \Gamma^*) Y_i^{} Y_j^* (Y_j^{} \Gamma^* - \Gamma Y_j^*) \\
	&\qquad + 2 \E \Gamma(I_p - Y_i^* Y_i^{}) Y_j^* (Y_j^{} \Gamma^* - \Gamma Y_j^*) + 2 \E (\Gamma Y_i^* - Y_i^{} \Gamma^*) Y_i (I_p - Y_j^* Y_j^{}) \Gamma^* \\
	&= 4 \ip{I_p - Y_i^* Y_i^{}}{ I_p - Y_j^* Y_j^{} } I_r + \ip{Y_i^* Y_i^{}}{Y_j^* Y_j^{}} I_r + \tr(Y_i^{} Y_j^*) Y_i^{} Y_j^*  \\
	&\qquad+ 2 \ip{I_p - Y_i^* Y_i^{}}{Y_j^* Y_j^{}} I_r  + 2 \ip{Y_i^* Y_i^{}}{I_p - Y_j^* Y_j^{}} I_r \\
	&= \parens*{ 4 (p - r) + \normF{Y_i^{} Y_j^*}^2 } I_r + \tr(Y_i^{} Y_j^*) Y_i^{} Y_j^*.
\end{align*}
We have used the fact that $\tr(A B^*)$ is always real when $A$ and $B$ are Hermitian.
Similar calculations as before, combined with \eqref{eq:imag_bound}, yield
\begin{align*}
	(4p - 2r) \sum L_{ij} \ip{Y_i}{Y_j}
	&\leq \sum L_{ij} ( r \normF{Y_i^{} Y_j^*}^2 + \tr^2(Y_i^{} Y_j^*) ) \\
	&\leq \sum L_{ij} \parens*{ r \normF{Y_i^{} Y_j^*}^2 + \frac{1}{2} \tr^2( Y_i^{} Y_j^* + Y_j^{} Y_i^* ) - r \normF{Y_i^{} Y_j^*}^2 } \\
	&= 4r \sum L_{ij} \ip{Y_i}{Y_j} - \frac{1}{2} \sum A_{ij} \normF{Y_i - Y_j}^4.
\end{align*}
Thus, if $G$ is connected, $2p \geq 3r$ implies $Y_1 = \cdots = Y_n$, or, equivalently, $Y Y^* = Z Z^*$.
The key benefit to this choice of $\Ydt$ is that we exactly cancel (rather than drop) the $\normF{Y_i^{} Y_j^*}^2$ terms that arise in~\eqref{eq:imag_bound}.

\newlength{\figwidth}
\newlength{\figheight}
\ifSIAM
\else
\setlength{\figwidth}{0.42\textwidth}
\setlength{\figheight}{0.30\textwidth}
\fi
\begin{figure}[t!]
	\centering
	\begin{minipage}[t]{0.48 \textwidth}
		\centering
		\begin{subfigure}[t][\figheight][t]{\figwidth}
			\centering
			\includegraphics[width=0.98\textwidth]{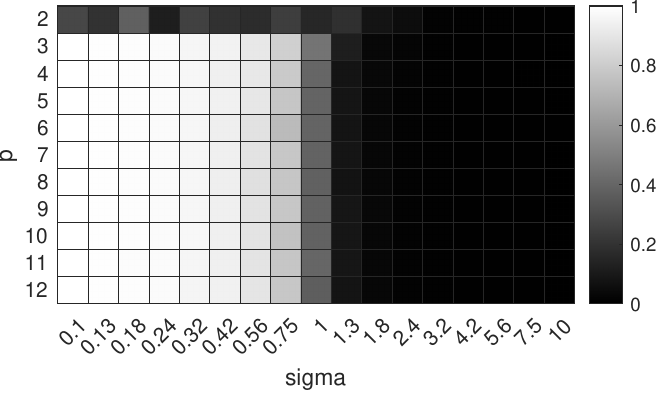}
			\subcaption{Normalized correlation $\normF{Z^\top Y}^2 / n^2 r$}
		\end{subfigure}\\
		\begin{subfigure}[t][\figheight][t]{\figwidth}
			\centering
			\includegraphics[width=0.98\textwidth]{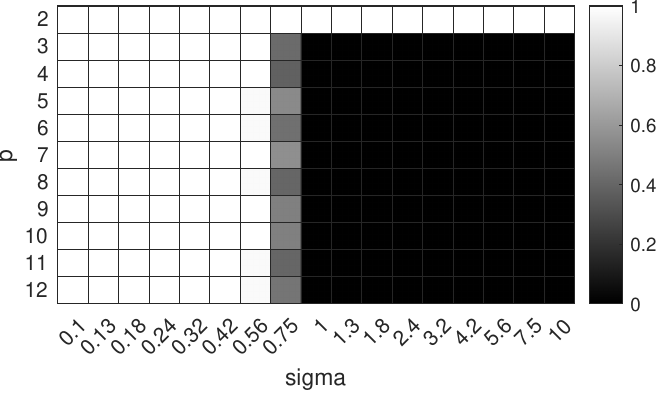}
			\subcaption{Solution is rank-$r$}
		\end{subfigure}\\
		\begin{subfigure}[t][\figheight][t]{\figwidth}
			\centering
			\includegraphics[width=0.98\textwidth]{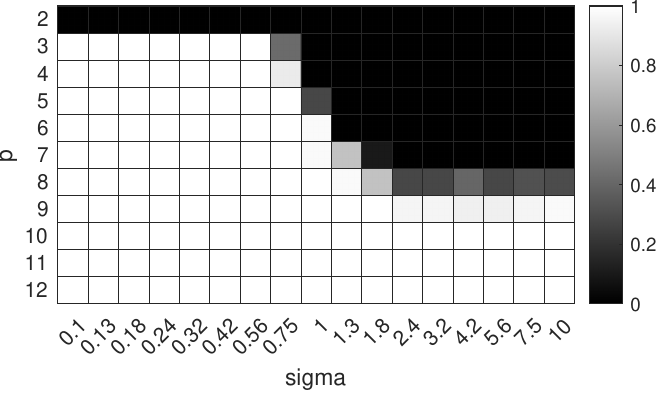}
			\subcaption{Solution is rank-deficient}
		\end{subfigure}\\
		\begin{subfigure}[t][\figheight][t]{\figwidth}
			\centering
			\includegraphics[width=0.98\textwidth]{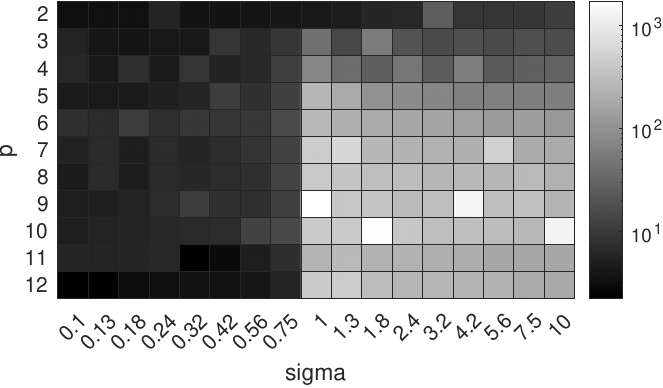}
			\subcaption{Solve time (s)}
		\end{subfigure}
		\caption{Experiments ($r = 2$) with a circulant graph on $n = 400$ vertices, each vertex having degree $10$. All values are the average of 50 experiments.}
		\label{fig:stats_circ}
	\end{minipage}\hfill	
	\begin{minipage}[t]{0.48 \textwidth}
		\centering
		\begin{subfigure}[t][\figheight][t]{\figwidth}
			\centering
			\includegraphics[width=0.98\textwidth]{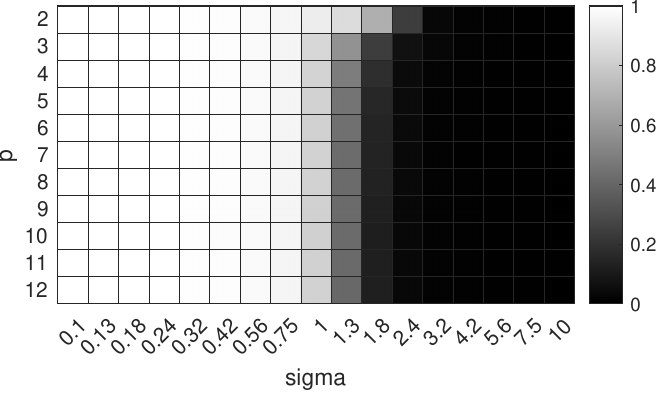}
			\subcaption{Normalized correlation $\normF{Z^\top Y}^2 / n^2 r$}
		\end{subfigure}\\
		\begin{subfigure}[t][\figheight][t]{\figwidth}
			\centering
			\includegraphics[width=0.98\textwidth]{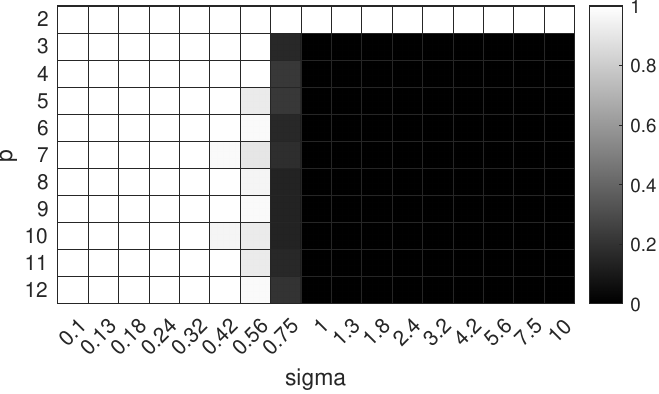}
			\subcaption{Solution is rank-$r$}
		\end{subfigure}\\
		\begin{subfigure}[t][\figheight][t]{\figwidth}
			\centering
			\includegraphics[width=0.98\textwidth]{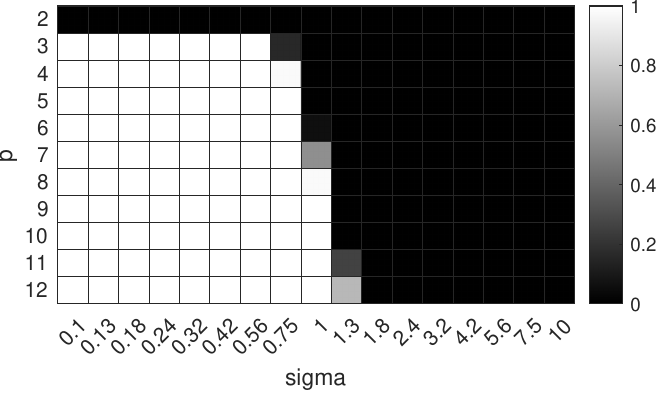}
			\subcaption{Solution is rank-deficient}
		\end{subfigure}\\
		\begin{subfigure}[t][\figheight][t]{\figwidth}
			\centering
			\includegraphics[width=0.98\textwidth]{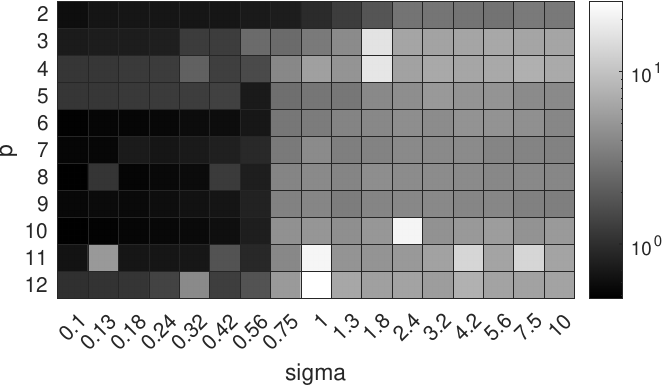}
			\subcaption{Solve time (s)}
		\end{subfigure}
		\caption{Experiments ($r = 2$) on a graph on $n = 400$ vertices chosen from an ER model with average degree $10$. All values are the average of 50 experiments.}
		\label{fig:stats_ER}
	\end{minipage}

\end{figure}

\begin{figure}[t]
	\centering
	\begin{subfigure}[t][\figheight][t]{\figwidth}
		\centering
		\includegraphics[width=0.98\textwidth]{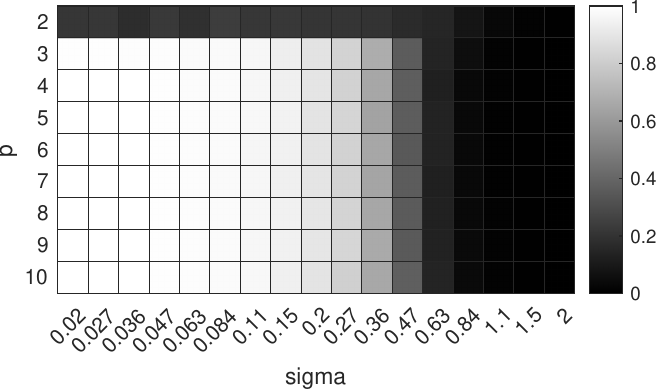}
		\subcaption{Normalized correlation $\normF{Z^\top Y}^2 / n^2 r$}
	\end{subfigure}
	\begin{subfigure}[t][\figheight][t]{\figwidth}
		\centering
		\includegraphics[width=0.98\textwidth]{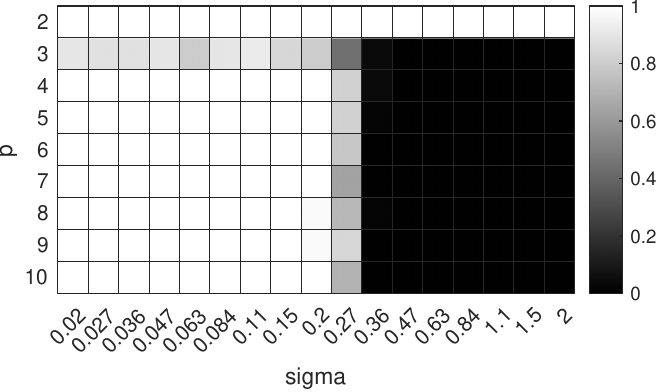}
		\subcaption{Solution is rank-$r$}
	\end{subfigure}\\
	\begin{subfigure}[t][\figheight][t]{\figwidth}
		\centering
		\includegraphics[width=0.98\textwidth]{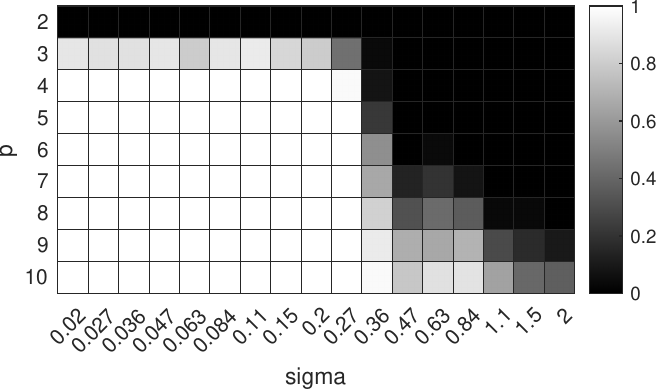}
		\subcaption{Solution is rank-deficient}
	\end{subfigure}
	\begin{subfigure}[t][\figheight][t]{\figwidth}
		\centering
		\includegraphics[width=0.98\textwidth]{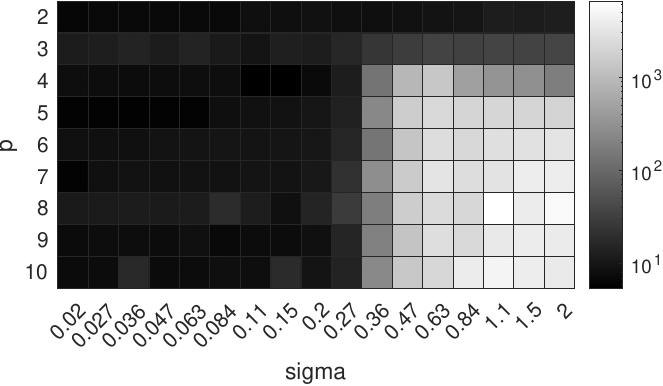}
		\subcaption{Solve time (s)}
	\end{subfigure}
	\caption{Experiments ($r = 2$) on the \texttt{fr079} SLAM dataset pose graph with $n = 989$ vertices of average degree $3.4$. All values are the average of 50 experiments.}
	\label{fig:stats_SLAM}
\end{figure}

\section{Simulations}
\label{sec:sims}
We implemented\footnote{Our code is available online at \texttt{\url{https://github.com/admcrae/sync_relaxed}}.} an algorithm for solving \eqref{eq:rankp_orig} and ran experiments on several graphs.
We used Matlab with the Manopt toolbox \cite{Boumal2014a} to optimize over a product of Stiefel manifolds with the default second-order trust-region algorithm.

We ran experiments on three graphs:
\begin{enumerate}
	\item A circulant graph on 400 vertices, each having degree 10 (results in \Cref{fig:stats_circ}).
	This graph topology is known to have spurious local minima (see, e.g., \cite{Townsend2020}).
	\item A single realization of the \ER{} random graph on 400 vertices, each vertex having expected degree 10 (\Cref{fig:stats_ER}).
	Because \ER{} graphs behave spectrally much like complete graphs, we expect them to be less prone to spurious local minima than the circulant graph (see \cite{Abdalla2022}).
	\item The pose graph from the Freiburg Building 079 (\texttt{fr079}) SLAM dataset in robotics,\footnote{Recorded by Cyrill Stachniss and available at \texttt{\url{https://www.ipb.uni-bonn.de/datasets/}}.} which has 989 vertices with average degree $3.4$ (\Cref{fig:stats_SLAM}).
\end{enumerate}
All experiments are in the real number case and use $r = 2$.
\Cref{thm:noiseless} predicts that, in the noiseless case, we should see benign nonconvexity for $p \geq 4$.
All experiments use random noise and initialization. In the case $p = r = 2$, the initial point is chosen to be in the same connected component as the ground truth.
The noise matrix $\Delta$ is chosen with i.i.d.\ $\normaldist(0, \sigma^2)$ entries in the nonzero blocks (with the constraint $\Delta = \Delta^\top$).
The reported $\rank(Y)$ is the number of singular values that are at least $10^{-3} \sqrt{n}$ (note that $Z$'s nonzero singular values are all $\sqrt{n}$).
The singular value tolerance did not qualitatively change the results (there was only a very slight effect at the phase transition).

The results are summarized in \Cref{fig:stats_circ,fig:stats_ER,fig:stats_SLAM}.
A few features worth highlighting are the following:
\begin{itemize}
	\item There is a clear phase transition as the noise standard deviation $\sigma$ increases;
	for each graph, the correlation performance, solution rank, and algorithm runtime all change dramatically in approximately the same place.
	
	\item The estimates cease to be rank-$r$ at only slightly lower noise levels than where the recovery performance noticeably begins to degrade.
	This suggests that, for these experiments, the prediction of \Cref{thm:noisy_landscape} (rank-$r$) is quite pessimistic compared to the known-to-be-optimal prediction of \Cref{thm:noisy_bound} (correlation).
	
	\item Except for the case $p = r = 2$ (see discussion below), the choice of $p$ has no noticeable effect on the correlation performance,
	even when the (local) optimum found is not rank-deficient.
	
	\item Other properties change markedly from one graph to another.
	\begin{itemize}
		\item \textbf{Noise tolerance:}
		The performance phase transition occurs soonest (i.e., for the smallest $\sigma$) for the large, sparse, and structured SLAM graph (\Cref{fig:stats_SLAM}).
		The highly structured circulant graph (\Cref{fig:stats_circ}) is next.
		The \ER{} graph (\Cref{fig:stats_ER}) can tolerate the most noise.
		
		\item \textbf{Runtime:}
		The runtimes follow the same trend as noise tolerance.
		The \ER{} graph (\Cref{fig:stats_ER}) has by far the shortest runtimes.
		The circulant graph (\Cref{fig:stats_circ}) has runtimes more than an order of magnitude higher,
		despite having the same size and (approximately) the same number of edges.
		The larger and sparser SLAM graph (\Cref{fig:stats_SLAM}) has by far the largest runtimes.
		\item \textbf{Solution rank:} For the circulant graph (\Cref{fig:stats_circ}), the solution rank is almost never larger than 8 even at very high noise levels. This is likely due to the simple graph structure.
		For the other graphs, the solution rank increases steadily for $\sigma$ near and above the phase transition.
		\item \textbf{Local optima for small $p$:} For the case $p = r = 2$ (which, in our experiments, is $\SOgrp(2)$ or angular synchronization),
		the \ER{} graph shows no problem with bad local optima for lower noise levels
		(curiously, the correlation performance is \emph{higher} than for larger $p$).
		This agrees with the Kuramoto oscillator results on \ER{} graphs such as \cite{Abdalla2022}.
		On the other hand, the algorithm is clearly finding bad local optima for the other, more structured graphs.
		Choosing $p = 3$ is adequate for the most part, though for the SLAM graph (\Cref{fig:stats_SLAM}),
		the algorithm also occasionally stops in a full-rank local optimum (but without any apparent loss of correlation performance).
	\end{itemize}

\end{itemize}

\section*{Acknowledgments}
We thank Pedro Abdalla, Afonso Bandeira, Johan Markdahl, David Rosen, and Alex Townsend for helpful conversations.

\ifSIAM
\else
\bibliographystyle{ieeetr}
\fi
\bibliography{./refs_siopt}
\end{document}